\documentclass[final,onefignum,onetabnum]{siamart190516}
\usepackage{lipsum}
\usepackage{amsfonts}
\usepackage{graphicx}
\usepackage{subfigure}
\usepackage{epstopdf}
\usepackage{algorithmic,version}
\ifpdf
  \DeclareGraphicsExtensions{.eps,.pdf,.png,.jpg}
\else
  \DeclareGraphicsExtensions{.eps}
\fi

\renewcommand\theequation{\thesection.\arabic{equation}}
\theoremstyle{plain}
\newtheorem{thm}{Theorem}[section]

\newtheorem{remark}{\textbf{Remark}}[section]
\newcommand{\bm}{\boldsymbol}

\newcommand{\Grad}[1]{\nabla #1}

\newcommand{\be}{\begin{equation}}
\newcommand{\ee}{\end{equation}}

\newcommand{\bse}{\begin{subequations}}
\newcommand{\ese}{\end{subequations}}
\def\benl{\begin{eqnarray*}}
\def\eenl{\end{eqnarray*}}

\def\be{\bm{e}}

\def\bx{\bm{x}}

\def\bmu1{\bm{\mu_1}}

\newcommand{\ben}{\begin{eqnarray}}
\newcommand{\een}{\end{eqnarray}}
\newcommand{\beq}{\begin{equation}}
\newcommand{\eeq}{\end{equation}}
\newcommand{\bea}{\begin{array}}
\newcommand{\eea}{\end{array}}
\newcommand{\bef}{\begin{figure}[H]}
\newcommand{\eef}{\end{figure}}

\newsiamremark{hypothesis}{Hypothesis}
\crefname{hypothesis}{Hypothesis}{Hypotheses}
\newsiamthm{claim}{Claim}

\title{Length preserving numerical schemes for Landau-Lifshitz equation based on Lagrange multiplier approaches}

\author{Qing Cheng\thanks{Department of  Mathematics, Purdue University  West Lafayette, IN 47907, USA 
  (\email{cheng573@purdue.edu)}. }
\and Jie Shen\thanks{Department of Mathematics, Purdue University, West Lafayette, IN 47907, USA  (shen7@purdue.edu)}. The work of J.S. is supported in part by NSF  DMS-1720442 and NFSC 11971407.}
\usepackage{amsopn}

\begin{document}
\bibliographystyle{plain}
\graphicspath{ {Figures/} }
\maketitle

\begin{abstract}
We develop in this paper two classes of length preserving schemes for  the Landau-Lifshitz equation based on  two different Lagrange multiplier approaches. In the first approach, the Lagrange multiplier $\lambda(\bx,t)$ equals to $|\nabla m(\bx,t)|^2$ at the continuous level, while in the second approach, the  Lagrange multiplier $\lambda(\bx,t)$ is introduced to enforce the length constraint at the discrete level and is identically zero at the continuous level. By using a predictor-corrector approach, we  construct efficient and robust length preserving higher-order schemes for the Landau-Lifshitz equation, with the computational cost dominated by the predictor step which is simply a semi-implicit scheme. Furthermore,  by introducing another space-independent Lagrange multiplier, we  construct energy dissipative, in addition to length preserving,  schemes  for the Landau-Lifshitz equation, at the expense of solving one nonlinear algebraic equation.
We present ample numerical experiments  to validate the stability and accuracy for the proposed  schemes, and also provide a performance comparison with some existing schemes.

\end{abstract}

% REQUIRED

% REQUIRED
%\begin{AMS}
 % 68Q25, 68R10, 68U05
%\end{AMS}


\section{Introduction}
We consider the Landau-Lifshitz equation (with $\beta\ne 0$, it is  often referred as the  Landau-Lifshitz-Gilbert equation) in the following form \cite{GuoH93,weinan2001numerical}:
\begin{equation}\label{Lf}
\begin{split}
&m_t = -\beta m\times \Delta m -\gamma m\times (m\times \Delta m)\quad  \mbox{in} \quad  \Omega,\\
&m(\bx,0)=m_0(\bx),
\end{split}
\end{equation}
with either homogeneous Neumann or periodic boundary condition. In the above, the unknown $m=(m_1,m_2,m_3)^t$ describes the magnetization in continuum ferromagnets,  $\Omega$ is an open bounded domain in $R^d$ $(d=1,2,3)$, $\beta$ is an exchange parameter, $\gamma>0$ is the Gilbert damping parameter,  and $m_0(\bx)$  with $|m_0(\bx)|\equiv 1$ is the initial condition. An important property of \eqref{Lf} is that the solution $m$ preserves pointwisely its magnitude. Indeed,
taking the dot product of \eqref{Lf} with $m$, we  derive   that vector field $m$ satisfies $\frac{d}{dt}|m(x)|=0$ for all $\bx\in \Omega$, so that solutions of  \eqref{Lf} satisfy an implicit constraint
\begin{equation}
	|m(\bx,t)|=1\quad\forall \bx\in \Omega,
\end{equation}	
i.e., $m$ is  length preserving.
Another important property is that \eqref{Lf} obeys an energy dissipation law. Taking the inner product of the first equation in \eqref{Lf} with $\Delta m$, we find that  \eqref{Lf}  satisfies the following energy dissipative law
\begin{equation}\label{energy}
	\frac{d}{dt}\int_{\Omega}\frac 12 |\Grad m|^2d\bx=-\gamma \|m\times \Delta m\|^2.
\end{equation}

The Landau-Lifshitz equation \eqref{Lf}, derived in \cite{landau1992theory} to describe the evolution of magnetization,  plays a very important role in understanding of non-equilibrium magnetism, and its accurate numerical simulation  has become an effective tool in understanding both the static and dynamics  in ferromagnetic materials \cite{kruzik2006recent,lakshmanan1984landau}.  Much effort has been developed in the past several decades to develop efficient  and accurate numerical methods for solving the Landau-Lifshitz equation \cite{gui2022convergence,an2021optimal,an2021analysis,gao2014optimal}.

To construct accurate and stable numerical schemes for the Landau-Lifshitz equation, it is of critical importance  to ensure that the key physical constraint $|m|=1$ is preserved at the discrete level. Existing numerical schemes for enforcing  the length constraint $|m|=1$  can be roughly classified into two categories: (i) Penalty approach \cite{pistella1999numerical,prohl2001computational}:   adding  a penalty term to approximate the constraint $|m|=1$. The penalty approach has been frequently used in numerical approximation of liquid crystal flows \cite{badia2011overview,badia2011finite,liu2000approximation,liu2002mixed}, but introduces additional numerical difficulties associated with the penalty parameters and does not enforce the length constraint exactly. (ii) Projection approach \cite{weinan2001numerical,bartels2008numerical,an2021optimal}: in the first-step, finding a approximate solution without enforcing the constraint, then preforming a simple projection in the second step to enforce the constraint.
Due to the simplicity of the projection approach, it has been frequently used, see, for instance. \cite{weinan2001numerical,wang2001gauss,alouges2006convergence,an2021optimal}.  Recently, a $L^2$-average orthogonal projection method is proposed in \cite{akrivis2021higher} where the length constraint is enforced in the sense of $L^2$-average. 
On the other hand, it appears to be difficult to construct higher-order robust schemes based on the projection approach. Furthermore, as our analysis indicates, the schemes based on the projection approach can not satisfy a discrete energy dissipation law.

In a sequence of recent work \cite{CS_CAMME22,CS_SINUM22}, we proposed efficient positivity/bound preserving schemes based on the Lagrange multiplier approach for a class of complex nonlinear systems.
The main purpose of this paper is to construct efficient length preserving schemes based on the Lagrange multiplier approach for the Landau-Lifshitz equation. More precisely, we introduce a new  Lagrange multiplier approach, which at its simplest form reduces to the simple projection approach proposed in \cite{weinan2001numerical}, and couple it with a predictor-corrector approach to construct efficient and robust length preserving higher-order schemes for the Landau-Lifshitz equation. Furthermore,  by introducing another space-independent Lagrange multiplier, we can also construct energy dissipative, in addition to length preserving,  schemes  for the Landau-Lifshitz equation, at the expense of solving one nonlinear algebraic equation. To the best of our knowledge, our schemes based on the predictor-corrector approach are the first  length preserving higher than second-order schemes for the Landau-Lifshitz equation, and our schemes with an additional space-independent Lagrange multiplier are the first  length preserving and energy dissipative schemes for the Landau-Lifshitz equation.

The paper is organized as follows. In Section $2$, we present  two different formulations of the  Landau-Lifshitz equation (for both the special case with $\beta=0$ and the general case with $\beta\ne 0$) by introducing  a Lagrange multiplier to enforce the length constraints.  In Section $3$, we construct a class of robust and accurate length preserving schemes for the special Landau-Lifshitz equation with $\beta=0$   based on the operator splitting and predictor-corrector approach.  In Section $4$, we modify  the length preserving schemes in Section 3 so that they also dissipate the energy. We then construct length preserving schemes and length preserving/energy decreasing schemes   for the general Landau-Lifshitz equation in Section 5. In Section $6$, we present ample numerical experiments  to validate the stability and accuracy for the proposed  schemes and provide a performance comparison with some existing schemes. We conclude with some remarks in the final section.

\section{Formulations of Landau-Lifshitz equation  with Lagrange multipliers}
It is pointed out in \cite{weinan2001numerical} that  descretizing directly \eqref{Lf:2b} without enforcing $|m|=1$ will lead to unstable numerical solutions. Hence, we shall consider expanded formulation of \eqref{Lf:2b}  by introducing a Lagrange multiplier to enforce the length constraint.

Since it is difficult to deal with $-m\times (m\times \Delta m)$ implicitly in a numerical scheme while an explicit treatment will lead to a severe time step constraint, we first
rewrite \eqref{Lf} as
\begin{equation}\label{Lf:2b}
	m_t = -\beta m\times \Delta m + \gamma(\Delta m + |\Grad m|^2 m).
\end{equation}
In the above, we used the fact that $-m\times (m\times \Delta m)=\Delta m + |\Grad m|^2 m$ if $|m|=1$, which can be derived from the identity
\begin{equation}
	a \times (b \times c) = (a\cdot c) b -(a\cdot b) c, \quad a,b,c \in R^3.
\end{equation}
Then, we introduce a Lagrange multiplier $\lambda(\bx,t)$ to enforce the length constraint and consider: 
\begin{equation}\label{Lf:2c}
	\begin{split}
		&m_t = -\beta m\times \Delta m + \gamma(\Delta m + |\Grad m|^2 m)+ \lambda(\bx,t)m,\\
		&|m(\bx,t)|=1.
	\end{split}
\end{equation}
Obviously with $\lambda(\bx,t)\equiv 0$, the above system is equivalent to \eqref{Lf:2b}. But we can now discretize \eqref{Lf:2c} directly with $|m(\bx,t)|=1$ explicitly enforced.

We shall also consider an alternative formulation  using a Lagrange multiplier as in \cite{weinan2001numerical}.
 Using the identity
\begin{equation}
	\frac 12	\Delta |m|^2= m\cdot\Delta m+|\nabla m|^2,
\end{equation}
we find that 
\begin{equation}\label{iden1}
	-m\cdot\Delta m=|\nabla m|^2 \quad \text{if }\; |m(x)|=1.
\end{equation}
Hence, we can rewrite \eqref{Lf:2b} with $|m|=1$ in the following equivalent form:
\begin{equation}\label{Lf:3b}
	\begin{split}
		&m_t = -\beta m\times \Delta m + \gamma(\Delta m +\lambda(\bx,t)m),\\
		&|m(\bx,t)|=1.
	\end{split}
\end{equation}
Indeed, multiplying the above by $m$ and using \eqref{iden1} and the fact that $|m(\bx,t)|=1$, we find that  $\lambda=|\Grad m|^2$. Hence, $\lambda$ can also be viewed as the  Lagrange multiplier  for the constraint $|m(\bx,t)|=1$.
%Hence, we shall discretize \eqref{Lf:3} with $|m|=1$ explicitly enforced instead of descretizing \eqref{Lf:2}. 

%Note that the Lagrange multiplier $\lambda(\bx,t)\equiv 0$ in \eqref{Lf:2c} is artificially introduced to enforce the implicit length constraint, while the Lagrange multiplier  in \eqref{Lf:3b} is $\lambda(\bx,t)=|\nabla m(\bx,t)|^2$.

A case of particular interest is when $\beta=0$. For the sake of simplicity, we also set $\gamma=1$. In this case, \eqref{Lf:2b} becomes
\begin{equation}\label{Lf:2d}
	m_t = \Delta m + |\Grad m|^2 m ;
	\end{equation}
while \eqref{Lf:2c} becomes
\begin{equation}\label{Lf:2}
\begin{split}
&	m_t = \Delta m + |\Grad m|^2 m +\lambda(\bx,t)m,\\
		&|m(\bx,t)|=1;
	\end{split}
\end{equation}
and \eqref{Lf:3b} becomes
\begin{equation}\label{Lf:3}
	\begin{split}
		&m_t =\Delta m +\lambda(\bx,t)m,\\
		&|m(\bx,t)|=1.
	\end{split}
\end{equation}
  The equation \eqref{Lf:2} or  \eqref{Lf:3} is also called the heat flow for harmonic maps, and has been extensively studied mathematically and numerically \cite{gui2022convergence,chen2021convergence,lin2008analysis}. 
  
  In the next section, we shall develop several numerical schemes for this special case based on both \eqref{Lf:2} and \eqref{Lf:3}. Extensions to the more general cases  based on \eqref{Lf:2c} and \eqref{Lf:3b} will be considered in Section 5. In the following, we shall refer to schemes based on \eqref{Lf:3}  and \eqref{Lf:3b}  as Type-I schemes, and those based on \eqref{Lf:2}  and  \eqref{Lf:2c} as Tyep-II schemes. Since the constructions of the Type-I and Type-II  schemes  follow essentially the same procedure, we shall present the Type-I schemes for  with sufficient details and explain briefly  how to construct Type-II  schemes.

\section{Length preserving time  discretization  schemes for the special Landau-Lifshitz equation \eqref{Lf:2d}}
In this section, we construct a class of efficient schemes for \eqref{Lf:3} and \eqref{Lf:2}. 
Note that $\lambda(\bx,t)$ (resp. $|m(\bx,t)|=1$) in  \eqref{Lf:3} and \eqref{Lf:2} plays a role similar to the pressure (resp. incompressibility constraint) in the Navier-Stokes equations. Hence, we can adopt the operator splitting and pressure-correction approaches developed for the  Navier-Stokes equations for solving    \eqref{Lf:3} and \eqref{Lf:2}.

\subsection{Type-I first-order operator-splitting scheme}
Similarly to the Chorin-Temam projection method for  the  Navier-Stokes equations \cite{Chor68,temam1968},  we  introduce a  first-order  operator-splitting scheme for Landau-Lifshitz equation \eqref{Lf:3}. Assuming  $m^n=(m_1^n,m_2^n,m_3^n)$ is known, we solve $\tilde m^{n+1}$ from
\begin{eqnarray}
\frac{\tilde{m}^{n+1}-m^n}{\delta t}=\Delta \tilde{m}^{n+1},\label{split:1}
\end{eqnarray}
and then we solve  $\lambda^{n+1}, m^{n+1}$ from
\begin{eqnarray}
&& \frac{m^{n+1}-\tilde{m}^{n+1}}{\delta t}=\lambda^{n+1}m^{n+1},\label{split:2}\\
&&|m^{n+1}|=1.\label{split:3}
\end{eqnarray}
\begin{thm} \label{proj}
The scheme \eqref{split:2}-\eqref{split:3}  admits two sets of solution, and the set of  solution consistent to \eqref{Lf:3} is given by
	\begin{equation}\label{project:2}
	\lambda^{n+1}=\frac{1-|\tilde{m}^{n+1}|}{\delta t},\quad	m^{n+1}=\frac{\tilde{m}^{n+1}}{|\tilde{m}^{n+1}|}.
	\end{equation}
\end{thm}

\begin{proof}
%\begin{remark}The scheme \eqref{split:1}-\eqref{split:3} can be comparable with the operator-splitting scheme for Naiver-Stokes equations \cite{} where equation \eqref{split:3} can be interpreted as  divergence free  condition in Naiver-Stokes equation.\end{remark}
 We rewrite \eqref{split:2} as 
\begin{equation}\label{step:1}
(1-\delta t \lambda^{n+1})m^{n+1}=\tilde{m}^{n+1}.
\end{equation}
Multiplying   \eqref{step:1} with itself on both sides, thanks to $|m^{n+1}|=1$,  we derive
\begin{equation}
(1-\delta t\lambda^{n+1})^2 = |\tilde m^{n+1}|^2.
\end{equation}
There are two roots for the above equation
\begin{equation}
\lambda^{n+1}=\frac{1-|\tilde{m}^{n+1}|}{\delta t} \quad \mbox{or} \quad \lambda^{n+1}=\frac{1+|\tilde{m}^{n+1}|}{\delta t}.
\end{equation}
Below we show that 
$\lambda^{n+1}=\frac{1-|\tilde{m}^{n+1}|}{\delta t}$
is the only  right solution. 

We derive from \eqref{split:1} that $\tilde{m}^{n+1}=(I-\delta t\Delta )^{-1}m^n\approx m^n+\delta t\Delta m^n$.   Hence, 
\begin{equation}
|\tilde{m}^{n+1}|^2 \approx 1+2\delta tm^n\cdot \Delta m^n+\delta t^2\delta t(m^n\cdot \Delta m^n)^2\approx (1+\delta t m^n\cdot \Delta m^n)^2\approx  (1-\delta t |\Grad m^n|^2)^2,
\end{equation}
where we used equality \eqref{iden1}.

Using the above equation, we derive
\begin{equation*}
\lambda^{n+1}=\frac{1-|\tilde{m}^{n+1}|}{\delta t}\approx\frac{1-(1-\delta t|\Grad m^n|^2)}{\delta t}\approx |\Grad m^n|^2,
\end{equation*}
which is consistent with  $\lambda(x,t)=-m(x,t)\cdot \Delta m(x,t)=|\Grad m|^2$ at the continuous level.
Plugging $\lambda^{n+1}=\frac{1-|\tilde{m}^{n+1}|}{\delta t}$ into equation \eqref{split:2}, we obtain $m^{n+1}=\frac{\tilde{m}^{n+1}}{|\tilde{m}^{n+1}|}$.
%\begin{equation}
%m^{n+1}=\frac{(I-\delta t\Delta )^{-1}m^n}{|(I-\delta t\Delta )^{-1}m^n|}.\end{equation}
\end{proof}

%		\section{Stability and error  analysis for the scheme \eqref{split:1}-\eqref{split:3} with $k=1$}
 The above result indicates that the scheme \eqref{split:1}-\eqref{split:3} is equivalent to the projection scheme in \cite{weinan2001numerical}. Hence, \eqref{split:1}-\eqref{split:3} is an alternative formulation of 	 the projection scheme, and it opens up a new avenue to develop higher-order version.
 
 	We recall that a stability and error  analysis in $L^\infty$-norm for the projection scheme  was carried out in \cite{weinan2001numerical}. Below,  we provide an alternative stability   analysis in $L^2$-norm for the scheme \eqref{split:1}-\eqref{split:3}. 
 		
 \begin{comment}
 		We recall the following discrete  Gronwall Lemma  \cite{shen1990long}:
 		\begin{lemma}\label{Gron2}
 			Let $a_n,\,b_n,\,c_n,$ and $d_n$ be four nonnegative sequences satisfying \begin{equation*}
 				a_m+\tau \sum_{n=1}^{m} b_n \le \tau \sum_{n=0}^{m-1}a_n d_n +\tau \sum_{n=0}^{m-1} c_n+ C, \,m \ge 1,
 			\end{equation*}
 			where $C$ and $\tau$ are two positive constants.
 			Then
 			\begin{equation*}
 				a_m+\tau \sum_{n=1}^{m} b_n \le \exp\big(\tau \sum_{n=0}^{m-1} d_n \big)\big(\tau \sum_{n=0}^{m-1}c_n+C \big),\,m \ge 1.
 			\end{equation*}
 		\end{lemma}
 \end{comment}		 
 		We first recall the following result from \cite{weinan2001numerical}.  
 		\begin{lemma}\label{lem:1}
 			Assume that 
 			\begin{equation}
 				(I-\alpha\delta t \Delta ) m = f \quad\text{in}\;\Omega 
 			\end{equation}
 			with the homogeneous Neumann boundary condition
 			or periodic boundary condition
 			where $m=(m_1,m_2,m_3)$, $f=(f_1,f_2,f_3)$ and $\alpha$ is any positive constant. Then
 			\begin{equation}
 				\max_{\forall x} |m| \leq \max_{\forall x} |f|.
 			\end{equation}
 		\end{lemma}
		\begin{comment}
 		\begin{proof}
 			For the readers' convenience, we provide a proof below. By direct calculation, we have 
 			\begin{equation}
 				\Delta |m|=\frac{1}{|m|}\Big[<m,\Delta m>+|\Grad m|^2-\frac{|<m,\Grad m>|^2}{|m|^2}\Big].
 			\end{equation}
 			Then 
 			\begin{equation}
 				\begin{split}
 					&(I-\alpha\delta t\Delta ) |m|=|m|-\alpha\delta t\Delta |m|\\&=
 					\frac{1}{|m|}<m,m-\alpha\delta t\Delta m>-\frac{\alpha\delta t}{|m|}\Big[|\Grad m|^2-\frac{|<m,\Grad m>|^2}{|m|^2}\Big]
 					\\&=\frac{1}{|m|}<m,f>-\frac{\alpha\delta t}{|m|}\Big[|\Grad m|^2-\frac{|<m,\Grad m>|^2}{|m|^2}\Big]
 					\\&\leq \frac{1}{|m|}<m,f> \leq |f|.
 				\end{split}
 			\end{equation}
 			Using the strong maximum principle for the following  equation
 			\begin{equation}
 				(I-\alpha\delta t\Delta ) |m| \leq |f|,
 			\end{equation}
 			we obtain the desired result.
 		\end{proof}
 	\end{comment}

 		\begin{lemma}\label{lem:2}
 			If the exact solution $m$ of Landau-Lifshitz equations poss enough regularity,  for Lagrange multiplier $\lambda^{n+1}$, we have the following bound
 			\begin{equation}\label{lbound}
 				0\leq \lambda^{n+1}(\bx) \leq C_0,
 			\end{equation}
 			for any $\bx$ in the domain $\Omega$, where $C_0$ is a constant independent of $\delta t$.
 		\end{lemma}
 		\begin{proof}
 			Since $|m^n|=1$, we derive immediately from  \eqref{split:1} and Lemma \ref{lem:1} that
 			\begin{equation}
 				|\tilde{m}^{n+1}| \leq 1.
 			\end{equation}
 			Then we derive from the above and  \eqref{project:2} that $\lambda^{n+1}\ge 0$. It is shown in (3.33) of  \cite{weinan2001numerical} that  there exists $C_0>0$ such that
 			\begin{equation*}
 				|\tilde m^{n+1}| \ge 1-C_0\delta t,
 			\end{equation*}
 			which,  together with  \eqref{project:2}, implies that  $\lambda^{n+1}(\bx) \leq C_0$. 
 		\end{proof}
 		
		We are now in position to prove the following stability results.
 		
 		\begin{theorem}\label{bound}
 		For the scheme \eqref{split:1}-\eqref{split:3}, we have  
 			\begin{equation}
 				\|\Grad m^k\|^2+2\delta t\sum\limits_{n=0}^{k-1}(\|\Grad \tilde{m}^{n+1}-\Grad m^n\|^2+\|\Grad m^{n+1}-\Grad \tilde{m}^{n+1}\|^2) \leq  \tilde{C}\|\Grad m^0\|^2,
 			\end{equation}
 			where $\tilde C$ is a positive constant depending on $C_0$ and $T$.
 		\end{theorem}
 		\begin{proof}
 			Taking the inner product of equation \eqref{split:1} with $-\frac{\tilde{m}^{n+1}-m^n}{\delta t}$, we obtain
 			\begin{equation}\label{split:stab:1}
 				-\|\frac{\tilde{m}^{n+1}-m^n}{\delta t}\|^2=\frac{1}{2\delta t}(\|\Grad \tilde{m}^{n+1}\|^2-\|\Grad m^n\|^2+\|\Grad \tilde{m}^{n+1}-\Grad m^n\|^2).
 			\end{equation}
 			Taking the inner product of equation \eqref{split:2} with $-\Delta m^{n+1}$, we obtain
 			\begin{equation}\label{split:stab:2}
 				\frac{1}{2\delta t}(\|\Grad m^{m+1}\|^2-\|\Grad \tilde{m}^{n+1}\|^2+\|\Grad m^{n+1}-\Grad \tilde{m}^{n+1}\|^2)=-(\lambda^{n+1}m^{n+1},\Delta m^{n+1}).
 			\end{equation}
 			Summing up equation \eqref{split:stab:1} and equation \eqref{split:stab:2}, we obtain
 			\begin{equation}\label{split:stab:3}
 				\begin{split}
 					&-\|\frac{\tilde{m}^{n+1}-m^n}{\delta t}\|^2-(\lambda^{n+1}m^{n+1},\Delta m^{n+1})\\&=\frac{1}{2\delta t}(\|\Grad m^{n+1}\|^2-\|\Grad m^n\|^2+\|\Grad \tilde{m}^{n+1}-\Grad m^n\|^2+\|\Grad m^{n+1}-\Grad \tilde{m}^{n+1}\|^2).
 				\end{split}
 			\end{equation}
 			Since $|m^{n+1}|=1$, we have $m^{n+1}\Delta m^{n+1}+|\Grad m^{n+1}|^2=0$, which implies
 			\begin{equation}
 				-(\lambda^{n+1}m^{n+1},\Delta m^{n+1})=(\lambda^{n+1},|\Grad m^{n+1}|^2).
 			\end{equation}
 			We can then rewrite  \eqref{split:stab:3} as
 			\begin{equation}\label{split:stab:4}
 				\begin{split}
 					&-\|\frac{\tilde{m}^{n+1}-m^n}{\delta t}\|^2=\frac{1}{2\delta t}(\|\Grad m^{n+1}\|^2-\|\Grad m^n\|^2\\&+\|\Grad \tilde{m}^{n+1}-\Grad m^n\|^2+\|\Grad m^{n+1}-\Grad \tilde{m}^{n+1}\|^2) -(\lambda^{n+1},|\Grad m^{n+1}|^2).
 				\end{split}
 			\end{equation}
 			By Lemma \ref{lem:2}, we have $0\leq \lambda^{n+1} \leq C_0$. Summing up equation \eqref{split:stab:4} for $n=1,2,\cdots,k-1=\frac{T}{\delta t}-1$, we obtain
 			\begin{equation}
 				\begin{split}
 					&\|\Grad m^k\|^2+2\delta t\sum\limits_{n=0}^{k-1}(\|\Grad \tilde{m}^{n+1}-\Grad m^n\|^2+\|\Grad m^{n+1}-\Grad \tilde{m}^{n+1}\|^2)
 					\\&\leq \|\Grad m^0\|^2 +2\delta t C_0 \sum\limits_{n=0}^{k-1}\|\Grad m^{n+1}\|^2.
 				\end{split}
 			\end{equation}
 			Applying a discrete Gronwall Lemma (cf. Lemma $2$ in \cite{shen1990long}) to the above, we obtain
 			\begin{equation}
 				\|\Grad m^k\|^2+2\delta t\sum\limits_{n=0}^{k-1}(\|\Grad \tilde{m}^{n+1}-\Grad m^n\|^2+\|\Grad m^{n+1}-\Grad \tilde{m}^{n+1}\|^2) \leq  \tilde{C}\|\Grad m^0\|^2.
 			\end{equation}
 			
 		\end{proof}

\subsection{Type-I higher-order predictor-corrector schemes}
Theorem \ref{proj} shows that the first-order operator splitting scheme \eqref{split:1}-\eqref{split:3} is equivalent to the projection scheme introduced   in \cite{weinan2001numerical}.  Hence,  it provides an alternative interpretation of the projection scheme. More importantly, it opens up  a new avenue to develop higher-order schemes for \eqref{Lf:2} through a predictor-corrector approach as we show below. Note that it is difficult to construct higher-order schemes based on the projection scheme, see  \cite{weinan2001numerical} for an attempt on constructing a second-order scheme.
However, using an idea similar to the pressure-correction scheme for the Navier-Stokes equations (see, for instance, \cite{guermond2006overview}),  we can construct higher-order schemes for \eqref{Lf:3} through a predictor-corrector approach as follows.

{\bf Step 1} (Predictor): solve $\tilde{m}^{n+1}$ from 
\begin{equation}\label{correction:1}
\frac{\alpha_k\tilde{m}^{n+1}-A_k(m^n)}{\delta t}=\Delta \tilde{m}^{n+1}+B_{k-1}(\lambda^n m^n),
\end{equation}

{\bf Step 2} (Corrector): solve $(m^{n+1},\lambda^{n+1})$ from
\begin{eqnarray}
&& \frac{\alpha_k (m^{n+1}-\tilde{m}^{n+1})}{\delta t}=\lambda^{n+1}m^{n+1}-B_{k-1}(\lambda^nm^n),\label{correction:2}\\
&&|m^{n+1}|=1,\label{correction:3}
\end{eqnarray}
where  $\alpha_k,$ and $A_k$ are determined from the $k$th-order Backward Difference-Formulas (BDF) \cite{akrivis2015stability,CS_CAMME22,CS_SINUM22}, and $B_{k-1}$ is determined from the Adams-Bashforth extrapolation:

\noindent {\bf First-order:}
\begin{equation}\label{eq:bdf1}
\alpha_1=1, \quad A_1(m^n)=m^n,\quad B_0(\lambda^n m^n)= 0;
\end{equation}
\noindent {\bf Second-order:}
\begin{equation}\label{eq:bdf2}
\alpha_2=\frac{3}{2}, \quad A_2(m^n)=2m^n-\frac{1}{2}m^{n-1},\quad B_1(\lambda^nm^n)=\lambda^n m^n;
\end{equation}
\noindent {\bf Third-order:}
\begin{equation}\label{eq:bdf3}
\alpha_3=\frac{11}{6}, \quad A_3(m^n)=3m^n-\frac{3}{2}m^{n-1}+\frac{1}{3}m^{n-2},\quad B_2(\lambda^nm^n)=2\lambda^nm^n-\lambda^{n-1}m^{n-1}.
\end{equation}
%\noindent {\bf Fourth-order:}
%\begin{equation}\label{eq:bdf4}\alpha_4=\frac{25}{12}, \; A_4(m^n)=4m^n-3m^{n-1}+\frac{4}{3}m^{n-2}-\frac{1}{4}m^{n-3},\; B_3(\lambda^nm^n)=3\lambda^nm^n-3\lambda^{n-1}m^{n-1}+\lambda^{n-2}m^{n-2}.\end{equation}
The formula for $k=4,5,6$ can be derived similarly with Taylor expansions.

We observe that the first step is a usual $k$th-order Implicit-Explicit (IMEX) scheme for the first equation in \eqref{Lf:3}, while the second step is $k$th-order correction to enforce $|m|=1$.
Obviously,  $\tilde m^{n+1}$ can be easily obtained from the first step.
Below we show how to efficiently solve $(m^{n+1},\lambda^{n+1})$  in the second step. 
We rewrite \eqref{correction:2} into the equivalent form
\begin{equation}\label{high:step:1}
\frac{\alpha_k \Big(m^{n+1}-\big(\tilde{m}^{n+1}-\frac{\delta t}{\alpha_k}B_{k-1}(\lambda^nm^n)\big)\Big)}{\delta t}=\lambda^{n+1}m^{n+1},
\end{equation}
and rearrange it into 
\begin{equation}\label{itself}
(\alpha_k-\delta t\lambda^{n+1})m^{n+1}=\alpha_k\tilde{m}^{n+1}-\delta tB_{k-1}(\lambda^nm^n).
\end{equation}
Multiplying \eqref{itself}  with  itself on both sides and using \eqref{correction:3}, we obtain
\begin{equation}
(\alpha_k-\delta t\lambda^{n+1})^2=|\alpha_k\tilde{m}^{n+1}-\delta tB_{k-1}(\lambda^nm^n)|^2.
\end{equation}
Then, similarly to the proof of Theorem \ref{proj}, we can establish the following result:
\begin{thm} \label{projh}
The scheme \eqref{correction:1}-\eqref{correction:3}  admits two sets of solution, and the set of  solution consistent to \eqref{Lf:3} is given by
\begin{equation}\label{sol:lambda}
\lambda^{n+1}=\frac{\alpha_k-|\alpha_k\tilde{m}^{n+1}-\delta tB_{k-1}(\lambda^nm^n)|}{\delta t},
\end{equation}
and 
\begin{equation}\label{sol:m}
m^{n+1}=\frac{\alpha_k\tilde{m}^{n+1}-\delta tB_{k-1}(\lambda^nm^n)}{|\alpha_k\tilde{m}^{n+1}-\delta tB_{k-1}(\lambda^nm^n)|}.
\end{equation}
\end{thm}

%\begin{equation}\tilde{m}^{n+1}=(I-\frac{\delta t}{\alpha_k}\Delta)^{-1}(\frac{1}{\alpha_k}A_k(m^n)+\frac{\delta t}{\alpha_k}B_{k-1}(\lambda^n m^n)).\end{equation}
%It is observed from \eqref{sol:m}, the second step actually represents  the  normalized process. 
In summary,  the scheme \eqref{correction:1}-\eqref{correction:3} can be implemented  as follows:
\begin{itemize}
\item Solve $\tilde{m}^{n+1}$ from \eqref{correction:1};
\item Find $m^{n+1}$ from projection step \eqref{sol:m};
\item Update $\lambda^{n+1} $ by \eqref{sol:lambda}.
\end{itemize}

\begin{remark}
We observe from \eqref{sol:m} that the scheme  \eqref{correction:1}-\eqref{correction:3} can be viewed as a modified projection method:

{\bf Step 1} (Predictor): solve $\tilde{m}^{n+1}$ from 
\begin{equation}\label{high:projection:1}
\frac{\alpha_k\tilde{m}^{n+1}-A_k(m^n)}{\delta t}=\Delta \tilde{m}^{n+1}+B_{k-1}(\lambda^n m^n),
\end{equation}
{\bf Step 2} (Projection): solve $m^{n+1}$ from 
\begin{equation}\label{high:projection:2}
m^{n+1}=\frac{\alpha_k\tilde{m}^{n+1}-\delta tB_{k-1}(\lambda^nm^n)}{|\alpha_k\tilde{m}^{n+1}-\delta tB_{k-1}(\lambda^nm^n)|},
\end{equation}
and update the  Lagrange multiplier $\lambda^{n+1}=\frac{\alpha_k-|\alpha_k\tilde{m}^{n+1}-\delta tB_{k-1}(\lambda^nm^n)|}{\delta t}$.

\end{remark}

\begin{remark}
We emphasize that one can construct other higher-order schemes using the predictor corrector approach. For example, 
a second-order scheme based on Crank-Nicolson is as follows:
Find $\tilde m^{n+1}$ from 
\begin{eqnarray}
 \frac{\tilde m^{n+1}-m^n}{\delta t}=\Delta \frac{\tilde m^{n+1}+\tilde m^n}{2}
+\lambda^n m^n;\label{cn:1}
\end{eqnarray}
and find $\lambda^{n+1},\, m^{n+1}$ from
\begin{eqnarray} 
&&\frac{m^{n+1}-\tilde{m}^{n+1}}{\delta t}=\frac{\lambda^{n+1}m^{n+1}-\lambda^nm^n}{2},\label{cn:2}\\
&&|m^{n+1}|=1.\label{cn:3}
\end{eqnarray}
It is easy to see that  \eqref{cn:2}-\eqref{cn:3} can be solved similarly as \eqref{correction:2}-\eqref{correction:3}. 

\end{remark}

\begin{remark}
The key to prove the uniform bound in Theorem \ref{bound} is the bound for $\lambda$ in Lemma \ref{lem:2}. On the other hand, we also observe from the bound \eqref{lbound} and \eqref{split:stab:4} that the scheme can not be energy dissipative.

It appears difficult to establish similar bound as  \eqref{lbound} for the above higher-order schemes, consequently, proving a uniform bound as in  Theorem \ref{bound}  is also illusive.
\end{remark}

\subsection{Type-II length preserving schemes}

Similarly to  \eqref{correction:1}-\eqref{correction:3},  we can construct  higher-order predictor-corrector schemes for  \eqref{Lf:2} as follows.

\noindent{\bf Step 1} (Predictor): solve $\tilde{m}^{n+1}$ from 
\begin{equation}\label{type2:correction:1}
\frac{\alpha_k\tilde{m}^{n+1}-A_k(m^n)}{\delta t}=\Delta \tilde{m}^{n+1}+B_{k}(|\Grad  m^{n}|^2  m^n)+B_{k-1}(\lambda^n m^n),
\end{equation}
{\bf Step 2} (Corretor): solve $(m^{n+1},\lambda^{n+1})$ from
\begin{eqnarray}
&& \frac{\alpha_k (m^{n+1}-\tilde{m}^{n+1})}{\delta t}=\lambda^{n+1}m^{n+1}-B_{k-1}(\lambda^nm^n),\label{type2:correction:2}\\
&&|m^{n+1}|=1,\label{type2:correction:3}
\end{eqnarray}
where  $\alpha_k,$ the operators $A_k$ and $B_{k-1}$ $(k=1,2,3)$ are given by \eqref{eq:bdf1}-\eqref{eq:bdf3}.

Similarly to  \eqref{correction:2}-\eqref{correction:3}, we can show that the above scheme admits two sets of solution and the one which is consistent  to  \eqref{Lf:2} is given by:
\begin{equation}\label{type2:sol:m}
\begin{split}
&\lambda^{n+1}=\frac{\alpha_k-|\alpha_k\tilde{m}^{n+1}-\delta tB_{k-1}(\lambda^nm^n)|}{\delta t},\\
&m^{n+1}=\frac{\alpha_k\tilde{m}^{n+1}-\delta tB_{k-1}(\lambda^nm^n)}{|\alpha_k\tilde{m}^{n+1}-\delta tB_{k-1}(\lambda^nm^n)|}.
\end{split}
\end{equation}
We observe that the only difference between the scheme \eqref{type2:correction:1}-\eqref{type2:correction:3} and 
the scheme  \eqref{correction:1}-\eqref{correction:3} is that one has to compute an  extra explicit term $B_{k}(|\Grad  m^{n}|^2  m^n)$ in  \eqref{type2:correction:1}.

\section{Energy decreasing and length preserving  schemes for \eqref{Lf:2d}}
The schemes presented in the last section preserve the length constraint $|m|=1$, but we are unable to show that they are  energy decreasing. In fact, we do not aware of any schemes for the Landau-Lifshitz equation \eqref{Lf:3} which are both   energy decreasing and length conserving.  In this section, we construct  schemes which are both energy decreasing schemes and length preserving for Landau-Lifshitz equation \eqref{Lf:3} by introducing an extra  Lagrange multiplier $\xi(t)$, which is independent of spatial variables, to enforce  energy dissipation. More precisely,  we consider the following expanded system for \eqref{Lf:3}:
\begin{equation}\label{en:Lf:3}
	\begin{split}
		&m_t =\Delta m +\lambda(\bx,t)m,\\
		& m=\frac{m+\xi(t)}{|m+\xi(t)|},\\
%		& \frac{d}{dt}E(u) =-\int_{\Omega} |m_t|^2 d\bx,
& \frac{d}{dt}E(u) =-\int_{\Omega} |m\times \Delta m|^2d\bx,
	\end{split}
\end{equation}
where $E(u)=\frac 12\int_{\Omega} |\Grad m|^2d\bx$. Note that the Lagrange multiplier $\xi(t)$ is introduced to  enforce energy dissipation. Obviously,  with $\xi(t)\equiv 0$, the above system reduces to    \eqref{Lf:3}.

%\subsection{Schemes which dissipate the original energy}
Assuming $\tilde m^n$, $m^n$ and  $\lambda^n$ are known, a Type-I first-order scheme for \eqref{en:Lf:3} is as follows:
 
 {\bf Step 1} (Predictor): solve $\tilde{m}^{n+1}$ from 
\begin{equation}\label{en:correction:1b}
\frac{\tilde{m}^{n+1}-m^n}{\delta t}=\Delta{\tilde m^{n+1}} +\lambda^n m^n,
\end{equation}

{\bf Step 2} (Corrector): solve $(\hat m^{n+1},\lambda^{n+1})$ from
\begin{eqnarray}
&& \frac{ \hat m^{n+1}-\tilde{m}^{n+1}}{\delta t}={\lambda^{n+1}\hat m^{n+1}}-\lambda^n m^n,\label{en:correction:2b}\\
&&  |\hat m^{n+1}|=1,\label{en:correction:3b}
\end{eqnarray}

{\bf Step 3} (Preserving energy dissipation): solve $(m^{n+1},\xi^{n+1})$ from
\begin{eqnarray}
&& m^{n+1}=\frac{\hat m^{n+1}+\xi^{n+1}}{|\hat m^{n+1}+\xi^{n+1}|},\label{en:correction:4b}\\
 %\frac{E^{n+1}-E^n}{\delta t}=-\|\frac{\hat m^{n+1}- \hat m^n}{\delta t}\|^2,
&& \frac{E^{n+1}-E^n}{\delta t}=-\|\hat m^{n+1}\times \Delta \hat m^{n+1}\|^2,\label{en:correction:5b}
\end{eqnarray}
where the energy  approximation $E^{n+1}$ is defined by
 \begin{equation}\label{en:deb}
\begin{split}
E^{n+1}=\frac 12\int_{\Omega} |\Grad m^{n+1}|^2d\bx.
\end{split}
\end{equation}
 Note that  the first-two  steps \eqref{en:correction:1b}-\eqref{en:correction:3b} can be solved the same way as   the scheme \eqref{split:1}-\eqref{split:3}. In particular, the consistent solution to the second step is:
 	\begin{equation}\label{project:2c}
	\lambda^{n+1}=\frac{1-|\tilde{m}^{n+1}-\delta t\lambda^n m^n|}{\delta t},\quad	\hat m^{n+1}=\frac{\tilde{m}^{n+1}}{|\tilde{m}^{n+1}|}.
	\end{equation}
 It remains to determine $\xi^{n+1}$ and $m^{n+1}$ from the  Step 3. Plugging \eqref{project:2c} into \eqref{en:correction:4b} and \eqref{en:correction:5b}, we  obtain a nonlinear 
algebraic equation for $\xi^{n+1}$:
\begin{equation}\label{nonalger:eq}
F(\xi^{n+1}):=E^{n+1} -E^n+ \delta t \|\hat m^{n+1}\times \Delta \hat m^{n+1}\|^2, %\|\frac{\hat{m}^{n+1}-m^n}{\delta t}\|^2=0,
\end{equation}
where $E^{n+1}$ is defined by
\begin{equation}
E^{n+1}=\frac 12\int_{\Omega}|\Grad \frac{\hat m^{n+1}+\xi^{n+1}}{|\hat m^{n+1}+\xi^{n+1}|} |^2 d\bx.
\end{equation}
To solving the above nonlinear algebraic equation, we can use either the Newton iteration or the following secant method:
\begin{equation}\label{newton}
\xi_{k+1}=\xi_k -\frac{F(\xi_k)(\xi_k-\xi_{k-1})}{F(\xi_k)-F(\xi_{k-1})}.
\end{equation}
Since $\xi^{n+1}$ is an approximation to zero,  we can choose  $\xi_1=0$ and $\xi_0=-O(\delta t^2)$. In all our experiments,  \eqref{newton} converges in a few iterations so that the cost is negligible.
Once we obtain $\xi^{n+1}$,  we can update $m^{n+1}$ by \eqref{en:correction:4b}.

Obviously, the scheme \eqref{en:correction:1b}-\eqref{en:correction:5b}  is first-order accurate, and it is unconditionally  energy stable in the sense that of \eqref{en:correction:5b}.

%However, it is not obvious how to construct higher-order  energy decreasing and length conserving schemes based on \eqref{en:Lf:3}. The main difficulty is how to approximate the energy dissipation law to higher-order. Since we can only use a one-step scheme for the energy law if we want to preserve  energy dissipation, it appears that the only choice for second-order is the following  scheme based on Crank-Nicolson as follows:

A second-order energy length preserving and energy decreasing scheme based on Crank-Nicolson is as follows:

{\bf Step 1} (Predictor): solve $\tilde{m}^{n+1}$ from 
\begin{equation}\label{en:correction:1}
\frac{\tilde{m}^{n+1}-m^n}{\delta t}=\Delta \frac{\tilde m^{n+1}+m^n}{2}+\lambda^n m^n,
\end{equation}

{\bf Step 2} (Corrector): solve $(\hat m^{n+1},\lambda^{n+1})$ from
\begin{eqnarray}
&& \frac{ \hat m^{n+1}-\tilde{m}^{n+1}}{\delta t}=\frac{\lambda^{n+1}\hat m^{n+1}-\lambda^nm^n}{2},\label{en:correction:2}\\
&&  |\hat m^{n+1}|=1,\label{en:correction:3}
\end{eqnarray}

{\bf Step 3} (Preserving energy dissipation): solve $(m^{n+1},\xi^{n+1})$ from
\begin{eqnarray}
&& m^{n+1}=\frac{\hat m^{n+1}+\xi^{n+1}}{|\hat m^{n+1}+\xi^{n+1}|},\label{en:correction:4}\\
%&& \frac{E^{n+1}-E^n}{\delta t}=-\|\frac{\hat m^{n+1}- \hat m^n}{\delta t}\|^2,
&& \frac{E^{n+1}-E^n}{\delta t}=-\|\frac{\hat m^{n+1}+m^n}{2}\times \Delta \frac{\hat m^{n+1}+m^n}{2}\|^2.\label{en:correction:5}
\end{eqnarray}
Note that the solution procedure for the above scheme is essentially the same as that of \eqref{en:correction:1b}-\eqref{en:correction:5b}. 

%The first two steps are obviously second-order approximation to \eqref{Lf:3}. But the third-step does not obviously lead to a second-order approximation for $m^{n+1}$ and $E^{n+1}$. Our numerical results in the next section indicate that the above scheme is much more accurate than first-order, but not quite second-second.

Similarly, we can construct Type-II  energy decreasing and length conserving schemes for \eqref{Lf:2}.

\section{Extension to the more general Landau-Lifshitz equation \eqref{Lf:3b} with $\beta\ne0$}
In this section, we extend our Lagrange multiplier approach to construct Type-I schemes for  the more general Landau-Lifshitz equation \eqref{Lf:3b}. Similar Type-II schemes can be constructed for \eqref{Lf:2c}.

If we treat the additional term $\beta m\times \Delta m $ totally explicitly, all the schemes that we constructed in the last two sections for \eqref{Lf:3} can be directly extended to \eqref{Lf:3b}. Indeed, the Type-I $k$-th order BDF schemes for \eqref{Lf:3b} can be constructed as follows:

{\bf Step 1} (Predictor): solve $\tilde{m}^{n+1}$ from 
\begin{equation}\label{LLG:correction:1}
\frac{\alpha_k\tilde{m}^{n+1}-A_k(m^n)}{\delta t}=\gamma(\Delta \tilde{m}^{n+1}+B_{k-1}(\lambda^n m^n))-\beta B_{k-1}(m^n\times \Delta m^n),
\end{equation}

{\bf Step 2} (Corrector): solve $(m^{n+1},\lambda^{n+1})$ from
\begin{eqnarray}
&& \frac{\alpha_k (m^{n+1}-\tilde{m}^{n+1})}{\delta t}=\gamma(\lambda^{n+1}m^{n+1}-B_{k-1}(\lambda^nm^n)),\label{LLG:correction:2}\\
&&|m^{n+1}|=1,\label{LLG:correction:3}
\end{eqnarray}
where  $\alpha_k,$ and $A_k$ are defined as before.

 However, explicit treatment of $\beta m\times \Delta m $  may lead to a severe time step constraint. Below, we shall construct an efficient length preserving scheme combining a stabilization technique \cite{SY2010} coupled with a Gauss-Seidel approach \cite{wang2001gauss,li2020two}.
%\subsection{Norm preserving schemes}

We first  rewrite \eqref{Lf:3b} in the following  form
\begin{equation}\label{LLG:lg}
	\begin{split}
		&m_t  -S(\Delta m-\Delta m)+ \beta m\times \Delta m = \gamma(\Delta m+\lambda(\bx, t) m),\\
		&|m|=1,
	\end{split}
\end{equation}
where $S>0$ is a stabilization  constant which will help to stabilize the time discretization.

Similarly to the scheme  \eqref{correction:1}-\eqref{correction:3}, we construct the following Type-I second-order  length preserving schemes for \eqref{LLG:lg}:

{\bf Step 1} (Gauss-Seidel predictor): solve $(\tilde{m}_1^{n+1}, \tilde m_2^{n+1},\tilde m_3^{n+1})$ from 
\begin{equation}\label{LLG:correction:1b}
	\begin{split}
		&\frac{\tilde{m}_1^{n+1}-m_1^n}{\delta t} -(S+\gamma )\Delta \tilde m^{n+1}_1=\gamma \lambda^n m_1^n-\beta(\tilde m^{n,\dagger}_2\Delta \tilde m^{n,\dagger}_3-\tilde m^{n,\dagger}_3\Delta \tilde m^{n,\dagger}_2)-S\Delta \tilde m^{n,\dagger}_1,\\
		&\frac{\tilde{m}_2^{n+1}-m_2^n}{\delta t}-(S+\gamma )\Delta \tilde m^{n+1}_2=\gamma \lambda^n m_2^n+\beta(\tilde m^{n+\frac 12}_1\Delta \tilde m^{n,\dagger}_3-\tilde m^{n,\dagger}_3\Delta \tilde m^{n+\frac 12}_1)-S\Delta \tilde m^{n,\dagger}_2,\\
		&\frac{\tilde{m}_3^{n+1}-m_3^n}{\delta t}-(S+\gamma )\Delta \tilde m^{n+1}_3=\gamma \lambda^n m_3^n-\beta(\tilde m^{n+\frac 12}_1\Delta \tilde m^{n+\frac 12}_2-\tilde m^{n+\frac 12}_2\Delta \tilde m^{n+\frac 12}_1)-S\Delta \tilde m^{n,\dagger}_3,
	\end{split}
\end{equation}
where $m_k^{n,\dagger}=\frac 32 m_k^n-\frac 12m_k^{n-1}$ and $\tilde m_k^{n+\frac 12}=\frac{\tilde m_k^{n+1}+m_k^n} {2}$ for $k=1,2,3$.

{\bf Step 2} (Corrector): solve $(m^{n+1},\lambda^{n+1})$ from
\begin{eqnarray}
	&& \frac{ m^{n+1}-\tilde{m}^{n+1}}{\delta t}=\gamma\frac{\lambda^{n+1}m^{n+1}-\lambda^nm^n}{2},\label{LLG:correction:2b}\\
	&&|m^{n+1}|=1.\label{LLG:correction:3b}
\end{eqnarray}
Note that in the Step 1, a Gauss-Seidel approach is used to deal with the term $\beta m\times \Delta m$. It is shown in \cite{wang2001gauss} that the Gauss-Seidel approach for $\beta m\times \Delta m$ can improve the stability compared with the totally explicit treatment while  only requiring to solve a sequence of constant coefficient elliptic problems. The solution procedure for the second step is the same as   in  \eqref{correction:2}-\eqref{correction:3}.

We can also add an additional step to the above scheme to preserve the energy dissipation as in the last section. For instance, a second-order length-preserving and energy stable scheme for \eqref{LLG:lg} can be constructed as follows:

{\bf Step 1} (Gauss-Seidel predictor): exactly the same as \eqref{LLG:correction:1}.

{\bf Step 2} (Corrector): solve $(\hat m^{n+1},\lambda^{n+1})$ from
\begin{eqnarray}
	&& \frac{\hat m^{n+1}-\tilde{m}^{n+1}}{\delta t}=\gamma\frac{\lambda^{n+1}\hat m^{n+1}-\lambda^nm^n}{2},\label{seidel:LLG:correction:2}\\
	&&|\hat m^{n+1}|=1.\label{seidel:LLG:correction:3}
\end{eqnarray}

{\bf Step 3} (Preserving energy dissipation): solve $(m^{n+1},\xi^{n+1})$ from
\begin{eqnarray}
&& m^{n+1}=\frac{\hat m^{n+1}+\xi^{n+1}}{|\hat m^{n+1}+\xi^{n+1}|},\label{seidel:LLG:correction:4}\\
&& \frac{E^{n+1}-E^n}{\delta t}=-\gamma \|\frac{\hat m^{n+1}+m^n}{2}\times \Delta \frac{\hat m^{n+1}+m^n}{2}\|^2.\label{seidel:LLG:correction:5}
\end{eqnarray}

Note that the second step is the same as above, so it can be solved in the same way. On the other hand, the last step is exactly the same as in the scheme \eqref{en:correction:1b}-\eqref{en:correction:5b}.

\begin{remark}
For  \eqref{LLG:correction:1}, we use a Gauss-Seidel approach which can increase allowable time steps for stable computation, see \cite{wang2001gauss}. We only need  to solve  three  Possion type equations from \eqref{LLG:correction:1}, the computational cost is  the same as a usual semi-implicit scheme.

\end{remark}

%{\bf Step 2} (Corrector): solve $(m_1^{n+1}, m_2^{n+1}, m_3^{n+1}, \lambda^{n+1})$ from  \eqref{correction:2}-\eqref{correction:3}.

\section{Numerical results}
In this section, we implement numerical experiments to validate the  convergence rate, accuracy, stability for type-I and type-II schemes we constructed above. Numerical results are shown for Landau-Lifshitz equation \eqref{Lf}. In space we consider periodic boundary condition and use Fourer-Spectral method \cite{shen2011spectral,canuto2012spectral}.

%	We need to show the following:
%\begin{itemize}	\item Convergence rate using initial condition (6.1) for equation (2.1) with (i) $\beta=0$; (ii) $\beta\ne 0$, $\gamma\ne0$; (iii) $\beta\ne 0$, $\gamma=0$.
%	\item Comparison using  initial condition (6.3): 
%\end{itemize}

\subsection{Convergence rate with a known exact solution}
We shall test the convergence rate for  the general Landau-Lifshitz equation \eqref{LLG:lg}  with an external force % which is defined by
%\begin{equation}
%m_t  = \gamma(\Delta m+\lambda(\bx, t) m) +\beta m\times \Delta m + f.
%\end{equation}
so that the exact solution is 
\begin{equation}
\begin{split}
&m^e_1(x,y,t)=\sin(t+x)\cos(t+y),\\
&m^e_2(x,y,t)=\cos(t+x)\cos(t+y),\\
&m^e_3(x,y,t)=\sin(t+y).
\end{split}
\end{equation}
We set $\Omega=[0,2\pi)^2$ with periodic boundary conditions and use the Fourier spectral method with $128\times 128$ modes for spatial approximation so that the spatial discretization error is negligible.  To test  convergence rates, we calculate the average $L^{\infty}$ errors  which  is defined by
\begin{equation}\label{A:error}
\|m-m^e\|_{L^{\infty}}=\frac{\|m_1^e-m_1\|_{L^{\infty}}+\|m_2^e-m_2\|_{L^\infty}+\|m_3^e-m_3\|_{L^\infty}}{3},
\end{equation}
where $m=(m_1,m_2,m_3)$ and $m^e=(m_1^e,m_2^e,m_3^e)$ are the numerical  solution and exact solution. 

We observe from Table.\;\ref{table1} that the expected convergence rates are obtained for schemes \eqref{LLG:correction:1}-\eqref{LLG:correction:3} with $k=1,2,3$ where the parameters are chosen to be $\gamma=\beta=1$ in \eqref{LLG:lg}. For the case $\gamma=1$ and $\beta=0$ in \eqref{LLG:lg}, the Crank-Nicolson scheme \eqref{cn:1}-\eqref{cn:3}  also achieves the second order convergence from Table.\;\ref{table1}.

\begin{table}[ht!]
\centering
\begin{tabular}{r||c|c|c|c|c|c|c|c}
\hline
$\delta t$      & {BDF1}  & Order          & {BDF2}  & Order &{BDF3} & Order & {CN} &Order \\ \hline
$1.6\times 10^{-3}$  & $3.02E(-5)$ & $-$   &$3.89E(-6)$  & $-$  & $6.97E(-10)$ & $-$    
&$2.03E(-6)$& $-$ \\\hline
$8\times 10^{-4}$  & $1.51E(-5)$ & $1.00$   &$9.72E(-7)$  & $2.00$  & $9.09E(-11)$& $2.94$  & $5.06E(-7)$ & $2.01$   \\\hline
$4\times 10^{-4}$  & $7.89E(-6)$ & $0.94$   &$2.43E(-7)$  & $2.00$  & $1.20E(-11)$& $2.92$   & $1.26E(-7)$ & $2.00$  \\\hline
$2\times 10^{-4}$  & $3.95E(-6)$ & $0.99$     &$6.07E(-8)$ & $2.00$ & $1.53E(-12)$ & $2.97$  & $3.15E(-8)$ & $2.00$ \\\hline
$1\times 10^{-4}$   & $1.97E(-6)$ & $1.00$   &$1.51E(-8)$ &$2.00$ & $2.18E(-13)$ & $2.81$   & $7.87E(-9)$ & $2.00$ \\\hline
$5\times 10^{-5}$  & $9.86E(-7)$ & $0.99$  &$3.79E(-9)$ &$1.99$ & $9.70E(-14)$ & $1.16$   & $1.96E(-9)$ & $2.01$ \\ \hline
\end{tabular}
\vskip 0.5cm
\caption{Accuracy test: The average $L^{\infty}$ error  error between $m=(m_1,m_2,m_3)$ and the exact solution $m_e$ at $t=0.01$ using  BDF$k$ schemes   for $k=1,2,3$ and Crank-Nicolson scheme \eqref{cn:1}-\eqref{cn:3}.}\label{table1}
\end{table}

\subsection{Convergence rate with a unknown exact solution}
Next, we  test the convergence rate for Crank-Nicolson scheme \eqref{cn:1}-\eqref{cn:3}, and 
energy decreasing schemes  \eqref{en:correction:1b}-\eqref{en:correction:5b} and \eqref{en:correction:1}-\eqref{en:correction:5} for the special case \eqref{Lf:2d}
with  the initial condition  
\begin{equation}
\begin{split}
&m_1(x,y,0)=\cos(x)\cos(y)\sin(0.1),\\
&m_2(x,y,0)=\cos(x)\cos(y)\cos(0.1),\\
&m_3(x,y,0)=\sqrt{1-\cos^2(x)\cos^2(y)},
\end{split}
\end{equation}
in the domain $[0,2\pi)^2$. We also  use $128$ Fourier modes  in each direction for spatial approximation. The exact solution is unknown so we compute a
 reference solution  by fourth-order Runge-Kutta method %using  the form
%\begin{equation}
%m_t=-\gamma m\times (m\times \Delta m),
%\end{equation}
with a small time step $\delta t=10^{-6}$.  
In  Table.\;\ref{table2},  the average $L^\infty$ errors between numerical solutions and the reference solution are shown. We observe that the Crank-Nicolson scheme \eqref{cn:1}-\eqref{cn:3} and the  first-order energy decreasing scheme \eqref{en:correction:1b}-\eqref{en:correction:5b} achieve second-order and first-order, respectively, while the  energy decreasing scheme  \eqref{en:correction:1}-\eqref{en:correction:5} leads to essentially   second-order convergence rate but is not as accurate as the Crank-Nicolson scheme \eqref{cn:1}-\eqref{cn:3}.

\begin{table}[ht!]
\centering
\begin{tabular}{r||c|c|c|c|c|c}
\hline
$\delta t$      & {CN}  & Order   & {\eqref{en:correction:1b}-\eqref{en:correction:5b}} &order       & { \eqref{en:correction:1}-\eqref{en:correction:5} }  & Order  \\ \hline
$8\times 10^{-4}$  & $1.44E(-3)$ & $-$  & $2.22E(-3)$ & $-$    &$1.49E(-3)$  & $-$       \\\hline
$4\times 10^{-4}$  & $1.01E(-5)$ & $-$  & $6.14E(-4)$ & $-$    &$3.21E(-4)$  & $2.21$  \\\hline
$2\times 10^{-4}$  & $2.58E(-6)$ & $1.97$  & $3.02E(-4)$ & $1.02$    &$1.16E(-4)$ & $1.47$ \\\hline
$1\times 10^{-4}$   & $6.50E(-7)$ & $1.99$  & $1.46 E(-4)$ & $1.05$   &$3.40E(-5)$ &$1.77$    \\\hline
$5\times 10^{-5}$  & $1.63E(-7)$ & $1.99$ & $6.76E(-5)$ & $1.11$   &$5.13E(-6)$ &$2.72$    \\ \hline
$2.5\times 10^{-5}$  & $4.12E(-8)$ & $1.98$  & $2.97E(-5)$ & $1.18$  &$1.26E(-6)$ &$2.19$    \\ \hline
\end{tabular}
\vskip 0.5cm
\caption{Accuracy test: The average  $L^{\infty}$ error  between $m=(m_1,m_2,m_3)$ and the reference solution at $t=0.01$ using Crank-Nicolson scheme \eqref{cn:1}-\eqref{cn:3}, energy decreasing scheme \eqref{en:correction:1b}-\eqref{en:correction:5b}  and \eqref{en:correction:1}-\eqref{en:correction:5} .}\label{table2}
\end{table}

\subsection{Comparisons between various second-order  schemes}
In this subsection, we compare accuracy and stability of various second-order schemes for the special Landau-Lifshitz equation \eqref{Lf:2d}. We set  $\Omega =[-\frac 12,\frac 12)^2$, and  test the benchmark problem in \cite{an2021optimal}  with   the initial condition 
\begin{equation}\label{non:ini}
m(\bx,0)=(m_1,m_2,m_3)=
\begin{cases}
(0,0,-1)^T & \text{for} \quad   |\bx|\ge \frac 12,\\
(\frac{2x_1A}{A^2+|\bx|^2},\frac{2x_2A}{A^2+|\bx|^2},\frac{A^2-|\bx|^2}{A^2+|\bx|^2})^T & \text{for} \quad   |\bx| \leq \frac 12,
\end{cases}
\end{equation}
where $A=(1-2|\bx|)^4$.

We consider  the following  second-order schemes:
\begin{itemize}
\begin{comment}
\item First-order backward Euler scheme:
\begin{equation}\label{scheme:1}
\frac{m^{n+1}-m^n}{\delta t}=\Delta m^{n+1} +|\Grad m^n|^2m^n.
\end{equation} 

\item First-order projection scheme in \cite{weinan2001numerical}:
\begin{equation}\label{scheme:2}
\begin{split}
&\frac{\tilde{m}^{n+1}-m^n}{\delta t}=\Delta \tilde{m}^{n+1},\\
& m^{n+1}=\frac{m^{n+1}}{|m^{n+1}|}.
\end{split}
\end{equation}
\item First-order BDF1 scheme with projection
\begin{equation}\label{BDF1:projection}
\begin{split}
&\frac{\tilde{m}^{n+1}-m^n}{\delta t}=\Delta \tilde{m}^{n+1}+|\Grad m^n|^2 m^n,\\
& m^{n+1}=\frac{m^{n+1}}{|m^{n+1}|}.
\end{split}
\end{equation}

\item Type-I first-order projection scheme:
\begin{equation}\label{scheme:3}
\begin{split}
&\frac{\tilde{m}^{n+1}-m^n}{\delta t}=\Delta \tilde{m}^{n+1}+\lambda^nm^n,\\
& \frac{m^{n+1}-\tilde{m}^{n+1}}{\delta t}=\lambda^{n+1}m^{n+1}-\lambda^n m^n,\\
&|m^{n+1}|=1.
\end{split}
\end{equation}

\item Type-II-first-order projection scheme:
\begin{equation}\label{scheme:4}
\begin{split}
&\frac{\tilde{m}^{n+1}-m^n}{\delta t}=\Delta \tilde{m}^{n+1}+|\Grad m^n|^2 m^n+\lambda^nm^n,\\
& \frac{m^{n+1}-\tilde{m}^{n+1}}{\delta t}=\lambda^{n+1}m^{n+1}-\lambda^n m^n,\\
&|m^{n+1}|=1.
\end{split}
\end{equation}
\end{comment}
\item  	The usual second-order semi-implicit scheme
\begin{equation}\label{scheme:5}
 \frac{ m^{n+1}-m^n}{\delta t}=\Delta \frac{ m^{n+1}+m^n}{2}
+\frac 32|\Grad m^n|^2m^n-\frac 12|\Grad m^{n-1}|^2m^{n-1}.
\end{equation}

\item The second-order scheme in \cite{weinan2001numerical}:
\begin{equation}\label{scheme:6}
\begin{split}
& \frac{\tilde m^{n+1}-m^n}{\delta t}=\Delta \frac{\tilde m^{n+1}+m^n}{2} +\delta t\Grad(|\Grad m^n|^2)\Grad m^n;\\
& m^{n+1}=\frac{\tilde m^{n+1}}{|\tilde m^{n+1}|}.
\end{split}
\end{equation}

\item Type-I Crank-Nicolson  predictor-corrector scheme:
\begin{equation}\label{scheme:7}
\begin{split}
& \frac{\tilde m^{n+1}-m^n}{\delta t}=\Delta \frac{\tilde m^{n+1}+m^n}{2}
+\lambda^n m^n;\\
&\frac{m^{n+1}-\tilde{m}^{n+1}}{\delta t}=\frac{\lambda^{n+1}m^{n+1}-\lambda^nm^n}{2};\\
&|m^{n+1}|=1.
\end{split}
\end{equation}

\item Type-II Crank-Nicolson  predictor-corrector  scheme:
\begin{equation}\label{scheme:8}
\begin{split}
& \frac{\tilde m^{n+1}-m^n}{\delta t}=\Delta \frac{\tilde m^{n+1}+m^n}{2}
+\frac 32|\Grad m^n|^2m^n-\frac 12|\Grad m^{n-1}|^2m^{n-1}+\lambda^n m^n;\\
&\frac{m^{n+1}-\tilde{m}^{n+1}}{\delta t}=\frac{\lambda^{n+1}m^{n+1}-\lambda^nm^n}{2};\\
&|m^{n+1}|=1.
\end{split}
\end{equation}

\item Type-II BDF$2$ predictor-corrector scheme:

\begin{equation}\label{scheme:9}
\begin{split}
& \frac{3\tilde m^{n+1}-4m^n+m^{n-1}}{2\delta t}=\Delta \tilde m^{n+1}
+|\Grad m^{\dagger,n}|^2m^{n,\dagger}+\lambda^n m^n;\\
&\frac{3m^{n+1}-3\tilde{m}^{n+1}}{2\delta t}=\lambda^{n+1}m^{n+1}-\lambda^nm^n;\\
&|m^{n+1}|=1,
\end{split}
\end{equation}
where $m^{\dagger,n}=2m^n-m^{n-1}$.
\end{itemize}

To compare accuracy of these schemes, we obtain a reference solution by using the fourth-order Runge-Kutta method 
with a very small time step $\delta t=10^{-6}$ in time and $64^2$ Fourier modes in space. Table.\;\ref{table4} shows the $L^{\infty}$ error for different schemes with  $\delta t=10^{-4}$. We observe that the simple semi-discrete scheme 
\eqref{scheme:5} without length preserving is unstable, and Schemes
 \eqref{scheme:7}, \eqref{scheme:8} and \eqref{scheme:9}  achieve similar second-order accuracy and the results are   much better than that of Scheme \eqref{scheme:6}.

   \begin{table}[ht!]
\centering
\begin{tabular}{rccccccc}
\hline
$T$      &  Scheme-$\eqref{scheme:5}$ & Scheme-$\eqref{scheme:6}$       & Scheme-$\eqref{scheme:7}$  & Scheme-$\eqref{scheme:8}$  &  Scheme-\eqref{scheme:9}\\ \hline
$0.01$  & $3.8198E(-07)$ & $2.0663E(-06)$   &$ 3.3198E(-09)$  & $ 9.3788E(-08)$      & $7.8620E(-07)$\\\hline
$0.02$  & $1.7265E(-04)$ & $6.4858E(-06)$   &$3.9061E(-09)$  & $6.1995E(-08)$      &$2.3073E(-07)$\\\hline
$0.04$  & NaN& $ 1.6174E(-05)$     &$ 7.6123E(-09)$ & $3.3997E(-08)$ 
& $4.9371E(-08)$ \\\hline
$0.06$  & - & $2.8101E(-05)$     &$4.0602E(-08)$ & $6.1291E(-08)$ 
& $1.7317E(-08)$ \\\hline
$0.08$   & -& $4.8164E(-05)$   &$ 1.7553E(-07)$ &$1.9378E(-07)$ & $9.2875E(-09)$   \\\hline
$0.10$  & - & $7.8866E(-05)$     &$5.7183E(-07)$ & $5.8832E(-07)$ 
& $1.3193E(-08)$\\\hline
$0.12$  & - & $1.2487E(-04)$  &$1.4952E(-06)$ &$ 1.5104E(-06)$ & $1.0365E(-07)$   \\ \hline
$0.16$  & - & $3.1270E(-04)$     &$7.3706E(-06)$ & $7.3549E(-06)$ & $3.6430E(-06)$ \\\hline
%\hline
%$CPU (seconds)$ & $-$ & $-$  &$-$&$-$  \\\hline
\end{tabular}
\vskip 0.5cm
\caption{Comparison of  Schemes \eqref{scheme:5}-\eqref{scheme:9}.}\label{table4}
\end{table}

\begin{figure}[htbp]
\centering
\subfigure[Scheme \eqref{scheme:7} and \eqref{scheme:8} with $\delta t=10^{-4},10^{-5}$]{
\includegraphics[width=0.45\textwidth,clip==]{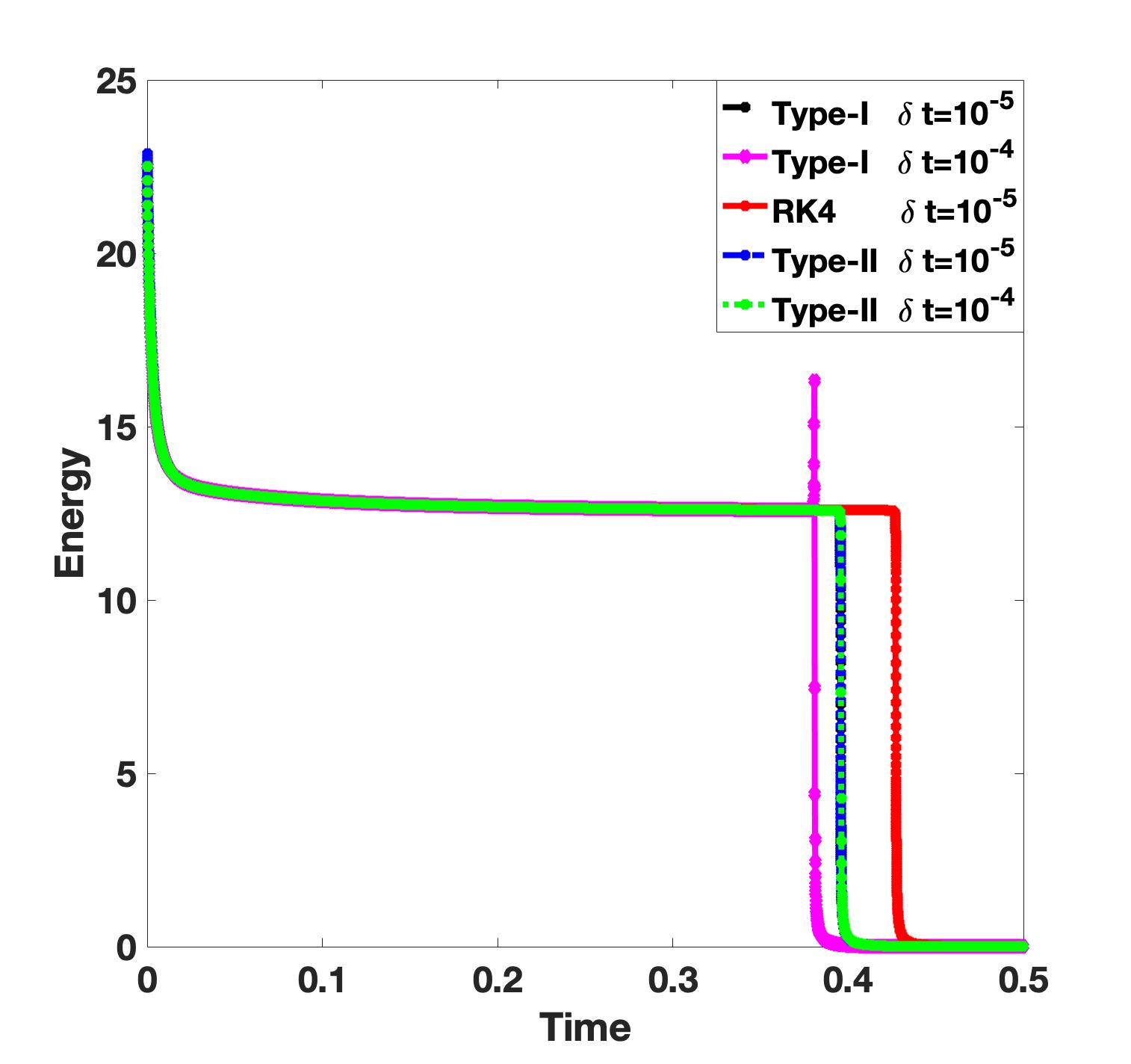}}
\subfigure[Type-II Scheme-\eqref{scheme:8} with $\delta t=10^{-3}$]{
\includegraphics[width=0.45\textwidth,clip==]{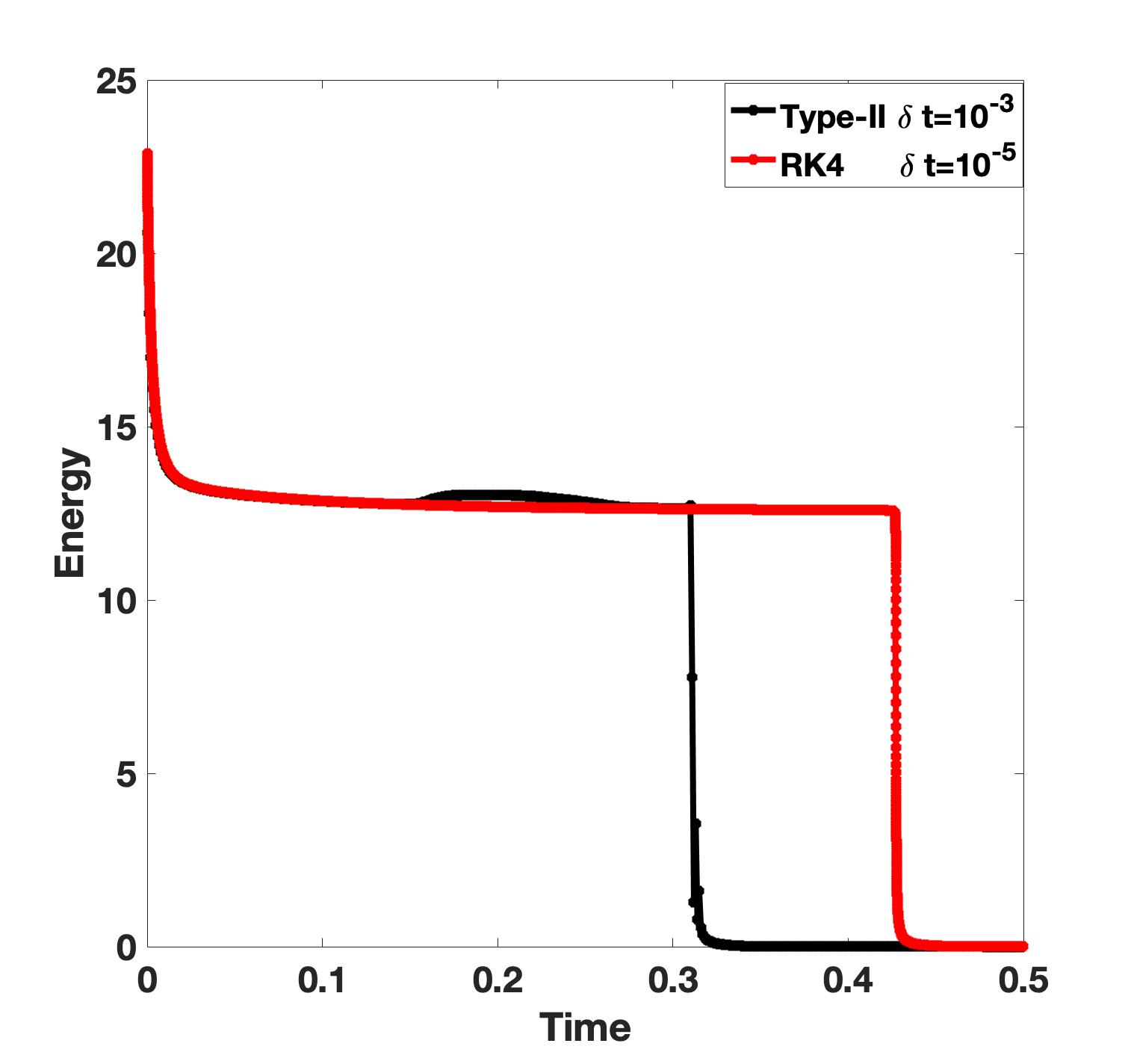}}
\subfigure[Type-I  Scheme-\eqref{scheme:7} with  $\delta t=10^{-3}$]{
\includegraphics[width=0.45\textwidth,clip==]{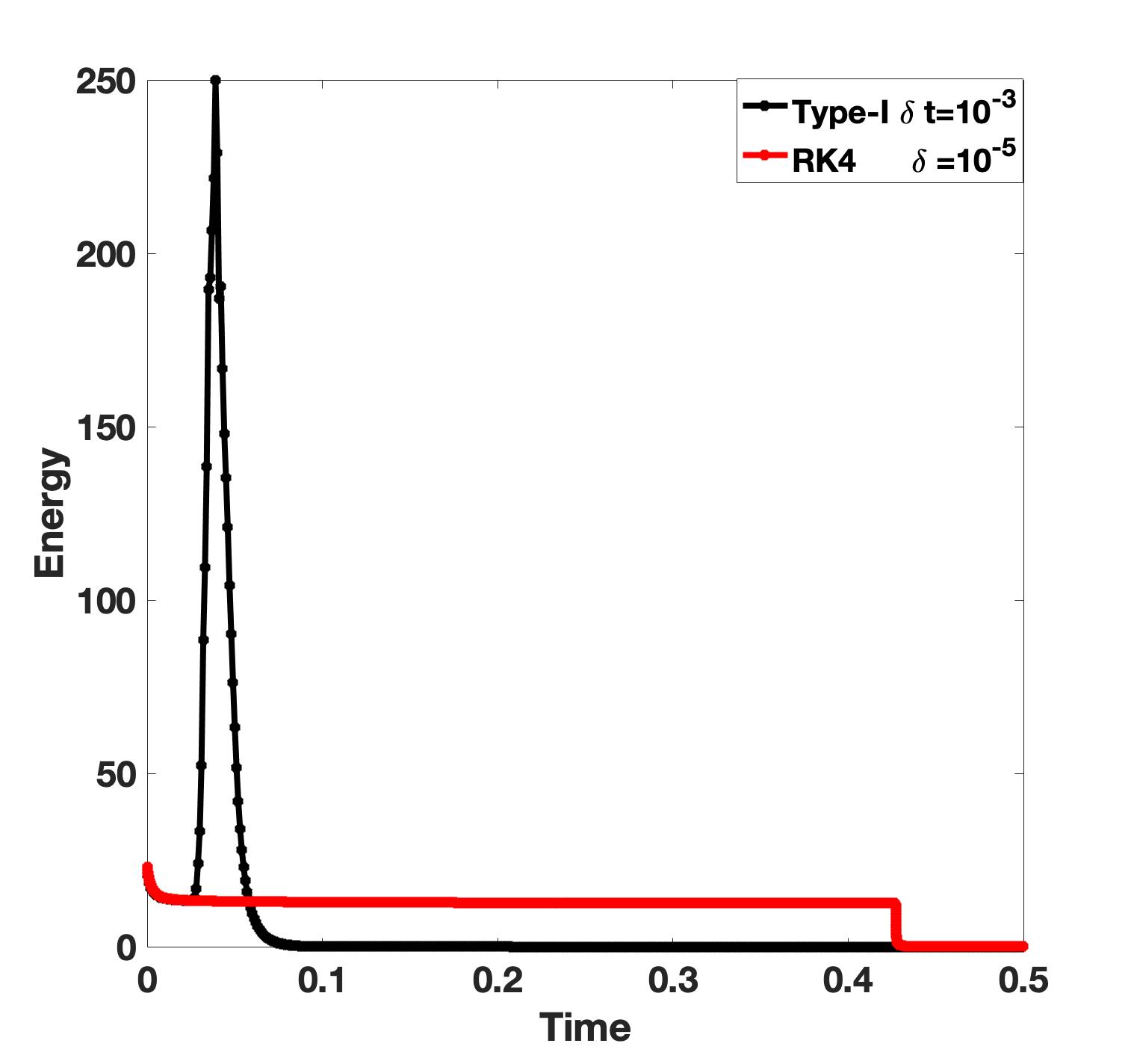}}
\caption{Energy comparison for Crank-Nicolson type-I  prediction-correction \eqref{scheme:7} and type-II prediction-correction scheme \eqref{scheme:8} using various time steps.}\label{LF_compare}
\end{figure}

\subsection{Comparison between Type-I and Type-II schemes}
In this subsection, we compare the stability and accuracy  between Type-I scheme \eqref{scheme:7}  and Type-II scheme \eqref{scheme:8} for Landau-Lifshitz equation \eqref{Lf:2d} with the initial condition given by \eqref{non:ini} in the domain $[-\frac 12,\frac 12)^2$. We also use Fourier-Spectral method in space with $N=64$ Fourier modes in each direction.

\begin{figure}[htbp]
\centering
\subfigure[$m$ at $t=0$.]{
\includegraphics[width=0.32\textwidth,clip==]{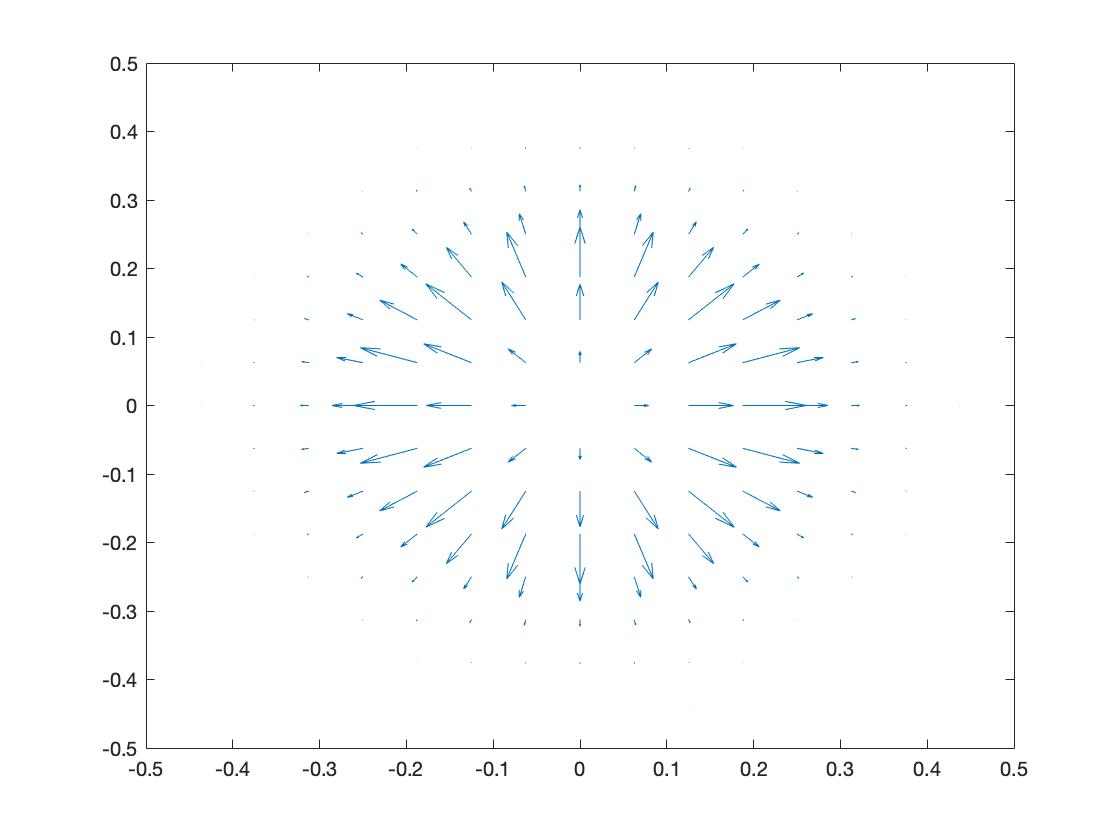}}
\subfigure[$m$ at $t=0.01$.]{
\includegraphics[width=0.32\textwidth,clip==]{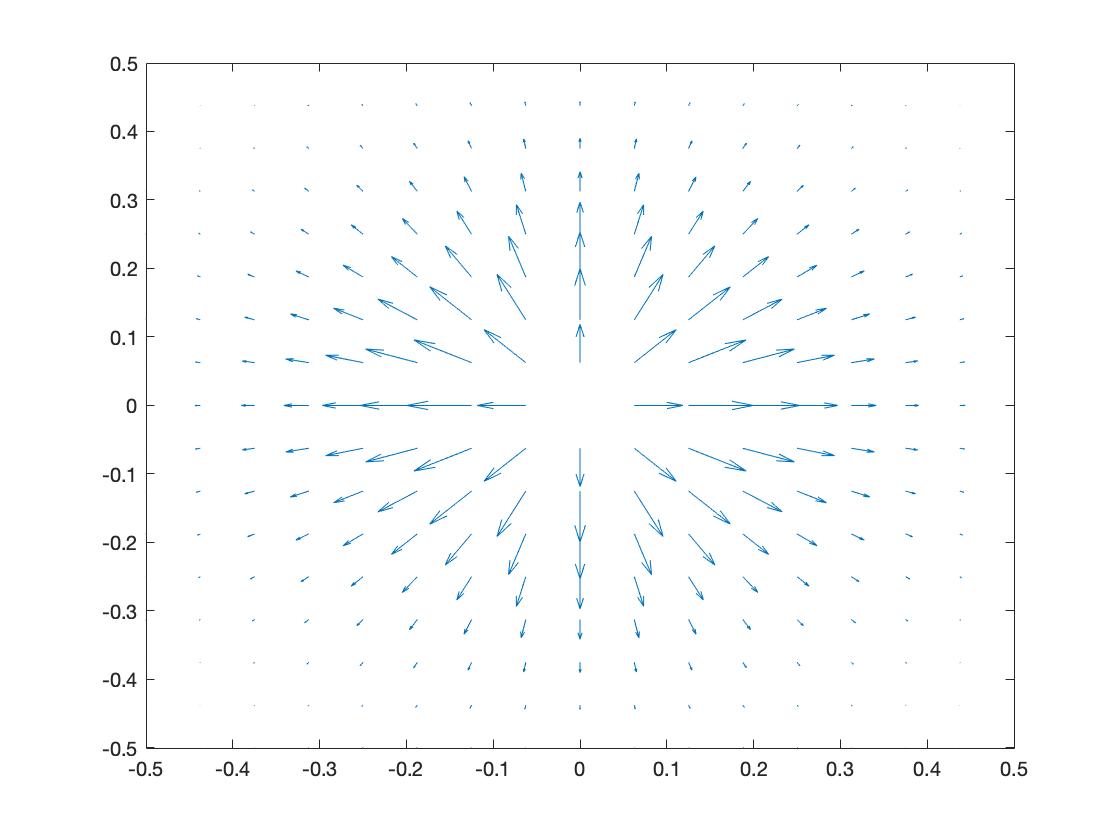}}
\subfigure[$m$ at $t=0.2$.]{
\includegraphics[width=0.32\textwidth,clip==]{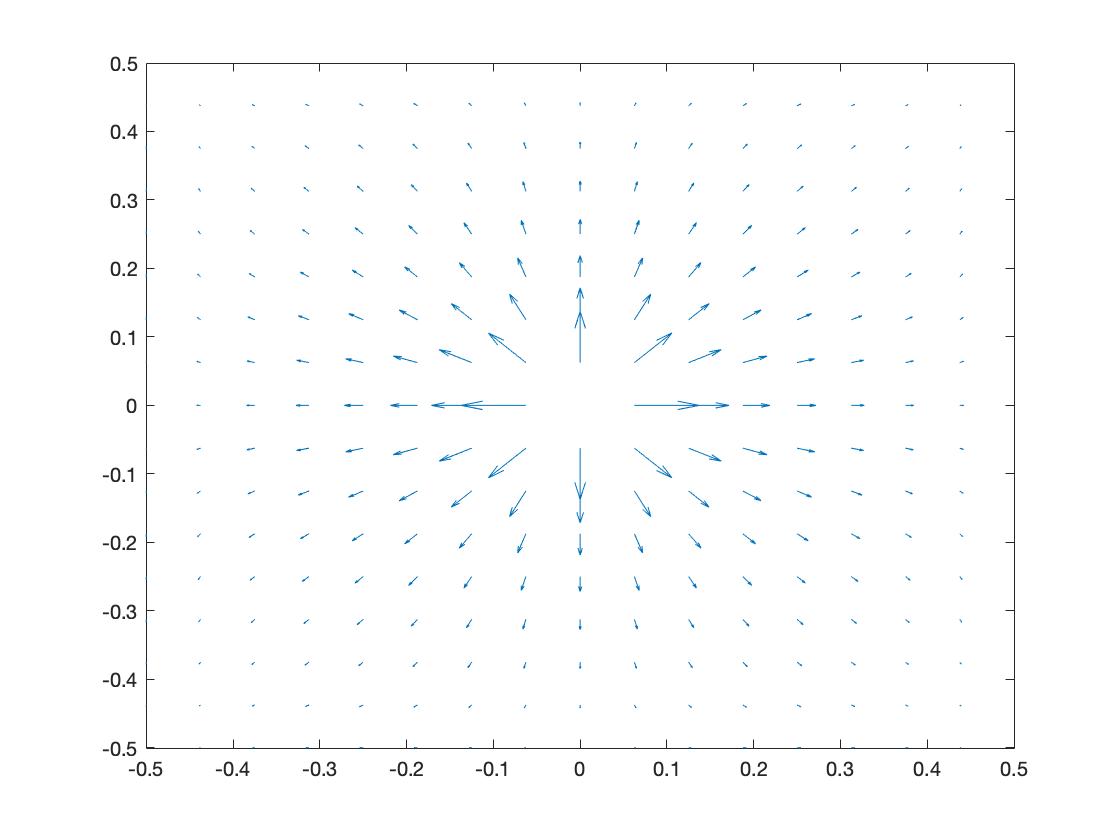}}
\subfigure[$m$ at $t=0.4$.]{
\includegraphics[width=0.32\textwidth,clip==]{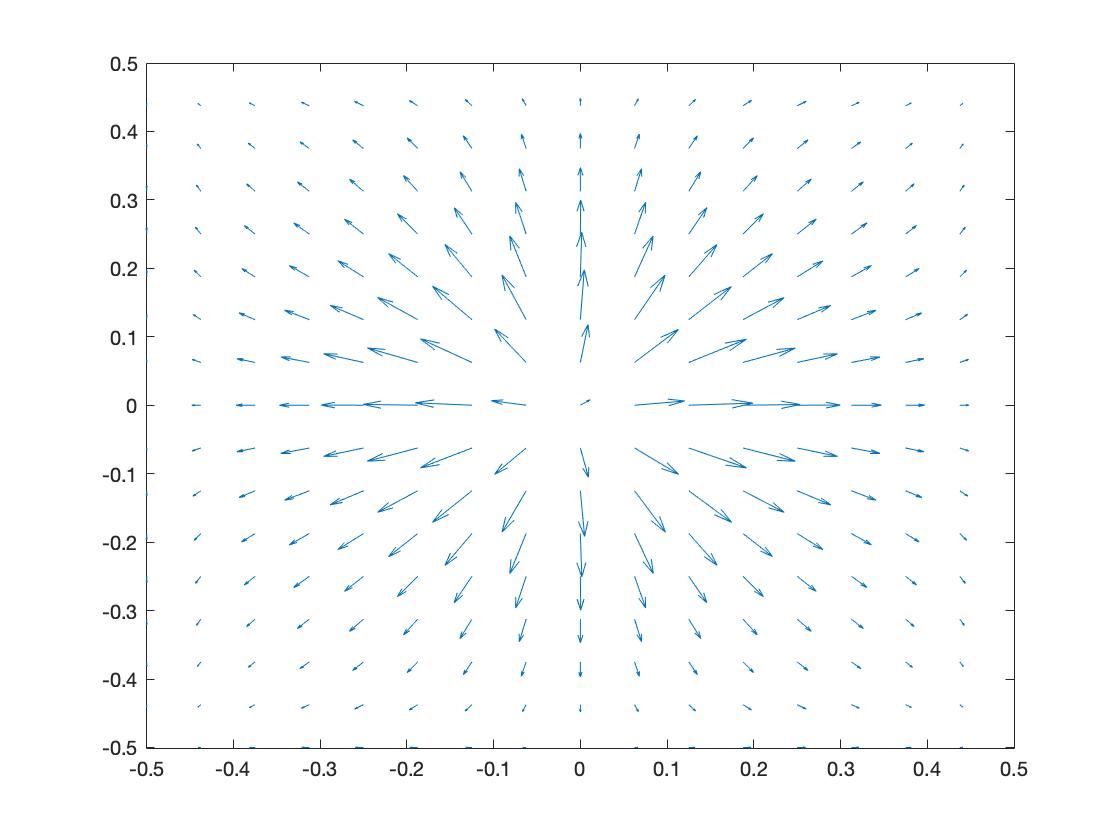}}
\subfigure[$m$ at $t=0.5$.]{
\includegraphics[width=0.32\textwidth,clip==]{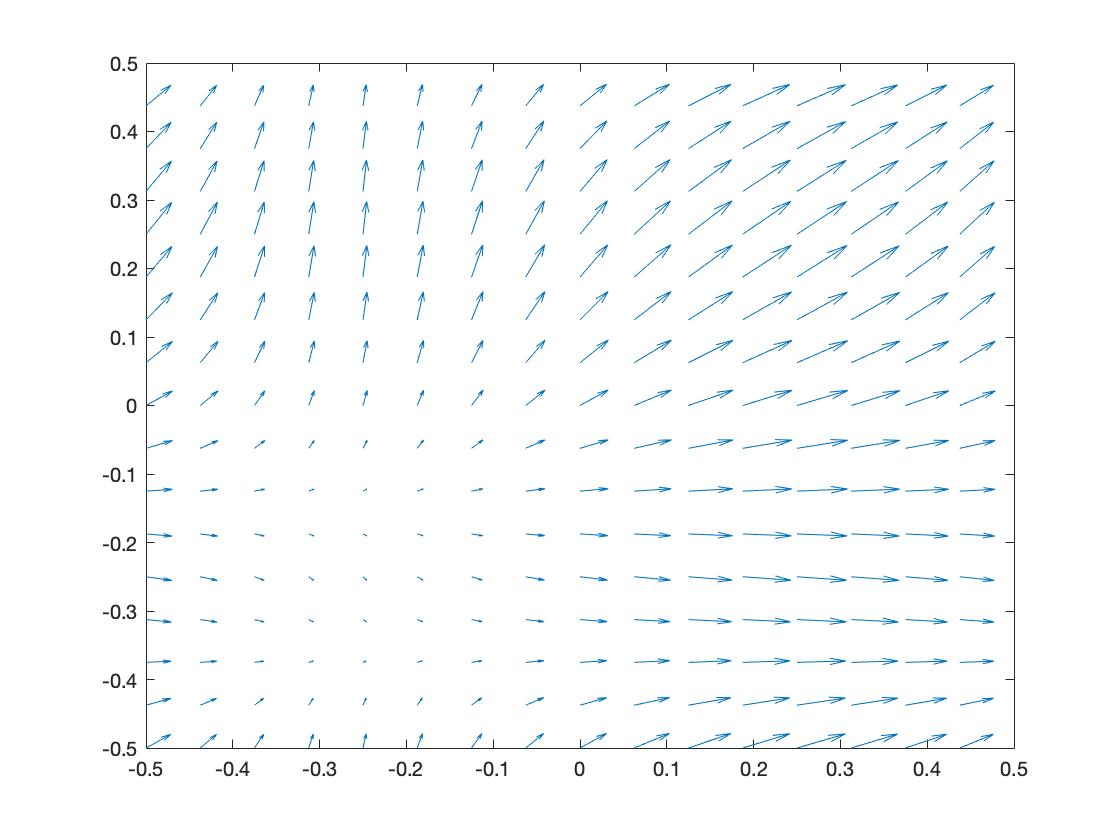}}
\subfigure[$m$ at $t=0.8$.]{
\includegraphics[width=0.32\textwidth,clip==]{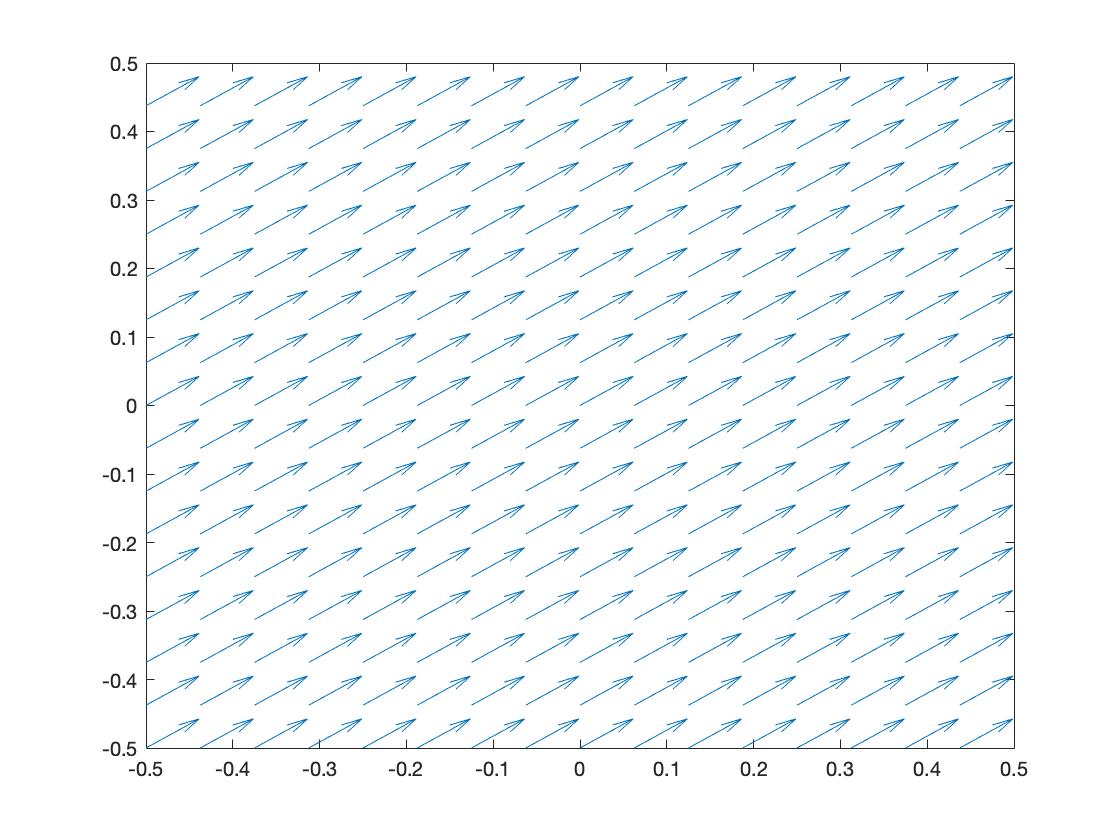}}
\caption{(a)-(f): Numerical solutions $m$ at $t=0, 0.01, 0.2, 0.4, 0.5, 0.8$ projected on $xy$-plane using the type-II Crank-Nicolson  scheme with $\delta t=10^{-4}$.}\label{LF_dynamic}
\end{figure}

In Fig.\;\ref{LF_compare}, we depict energy curves for numerical solutions  computed by  Scheme \eqref{scheme:7} and Scheme \eqref{scheme:8} with different time steps. In  Fig.\;\ref{LF_compare}.(a), we observe that both  schemes are stable  and can simulate  accurate dynamics of  Landau-Lifshitz equation \eqref{Lf:2d} with time steps $\delta t=10^{-4}, 10^{-5}$.  In Fig.\;\ref{LF_compare}.(b) and Fig.\;\ref{LF_compare}.(c), it is observed that the Type-II scheme \eqref{scheme:8} can produce more accurate  energy curve  than the  Type-I scheme  \eqref{scheme:7} at the larger time step $\delta t=10^{-3}$.  Numerical solutions $m$ at  $t=0, 0.01, 0.2, 0.4, 0.5, 0.8$ projected on $xy$-plane computed by the type-II Crank-Nicolson scheme are shown in Fig.\;\ref{LF_dynamic}. 

\subsection{Adaptive time stepping schemes}

A main advantage of unconditionally energy stable schemes is that one can adaptively choose time steps based on the accuracy only. Below, we present an adaptive algorithm which can be used with the unconditionally energy stable schemes introduced in previous section.

\medskip
\noindent{\bf Algorithm for adaptive time stepping:}\\
\rule[4pt]{14.3cm}{0.05em}\\
\textbf{Given}  Solutions at time steps $n$ and $n-1$; parameters $tol$, and the preassigned minimum and maximin allowable time steps $\delta t_{min}$ and $\delta t_{max}$.
\begin{description}
\item[Step 1] Compute $\hat m^{n+1}$  by Crank-Nicolson  scheme \eqref{en:correction:1}-\eqref{en:correction:3}  with $\delta t_n$.
\item[Step 2] Compute $(\xi^{n+1},m^{n+1})$  from \eqref{en:correction:4}-\eqref{en:correction:5}.
\item[Step 3] \textbf{if} $e^{n+1}=|\xi^{n+1}| > tol$, \textbf{then}\\
Recalculate time step $\delta t_n \leftarrow \max\{\delta t_{min},\min\{A_{dp}(e^{n+1},\delta t^n),\delta t_{max}\}\}$.
\item[Step 4] \textbf{goto} Step 1
\item[Step 5] \textbf{else}\\
Update time step  $\delta t_{n+1}\leftarrow \max\{\delta t_{min},\min\{A_{dp}(e^{n+1},\delta t^n),\delta t_{max}\}\}$,
\item[Step 6] \textbf{endif}
\end{description}
\rule[12pt]{14.3cm}{0.05em}\\
One simple but effective strategy is to update the time step size by using the formula \cite{cheng2018multiple},
\begin{equation}
A_{dp}(e,\delta t)=\rho(\frac{tol}{e})^{\frac 12}\delta t.
\end{equation}
We choose the constant $\rho=0.95$ for  numerical simulations in this subsection.
We implemented the above adaptive stepping based on  the energy decreasing scheme \eqref{en:correction:1}-\eqref{en:correction:5} for \eqref{Lf:2d}   and the energy decreasing scheme   \eqref{seidel:LLG:correction:2}-\eqref{seidel:LLG:correction:5} for \eqref{Lf}.

\begin{figure}[htbp]
\centering
\subfigure[evolution of $\xi^{n+1}$]{
\includegraphics[width=0.4\textwidth,clip==]{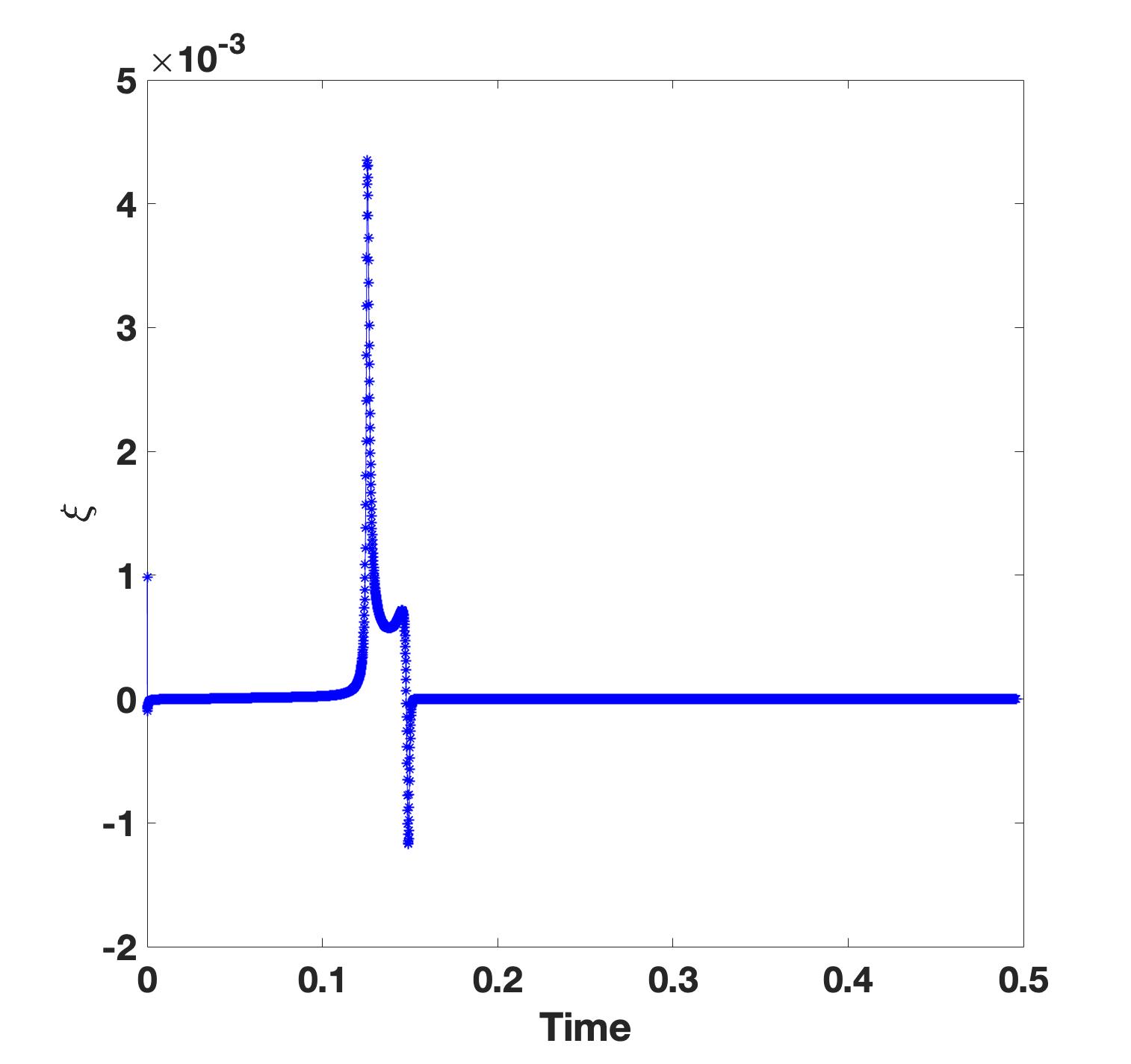}}
\subfigure[evolution of energy curves]{
\includegraphics[width=0.4\textwidth,clip==]{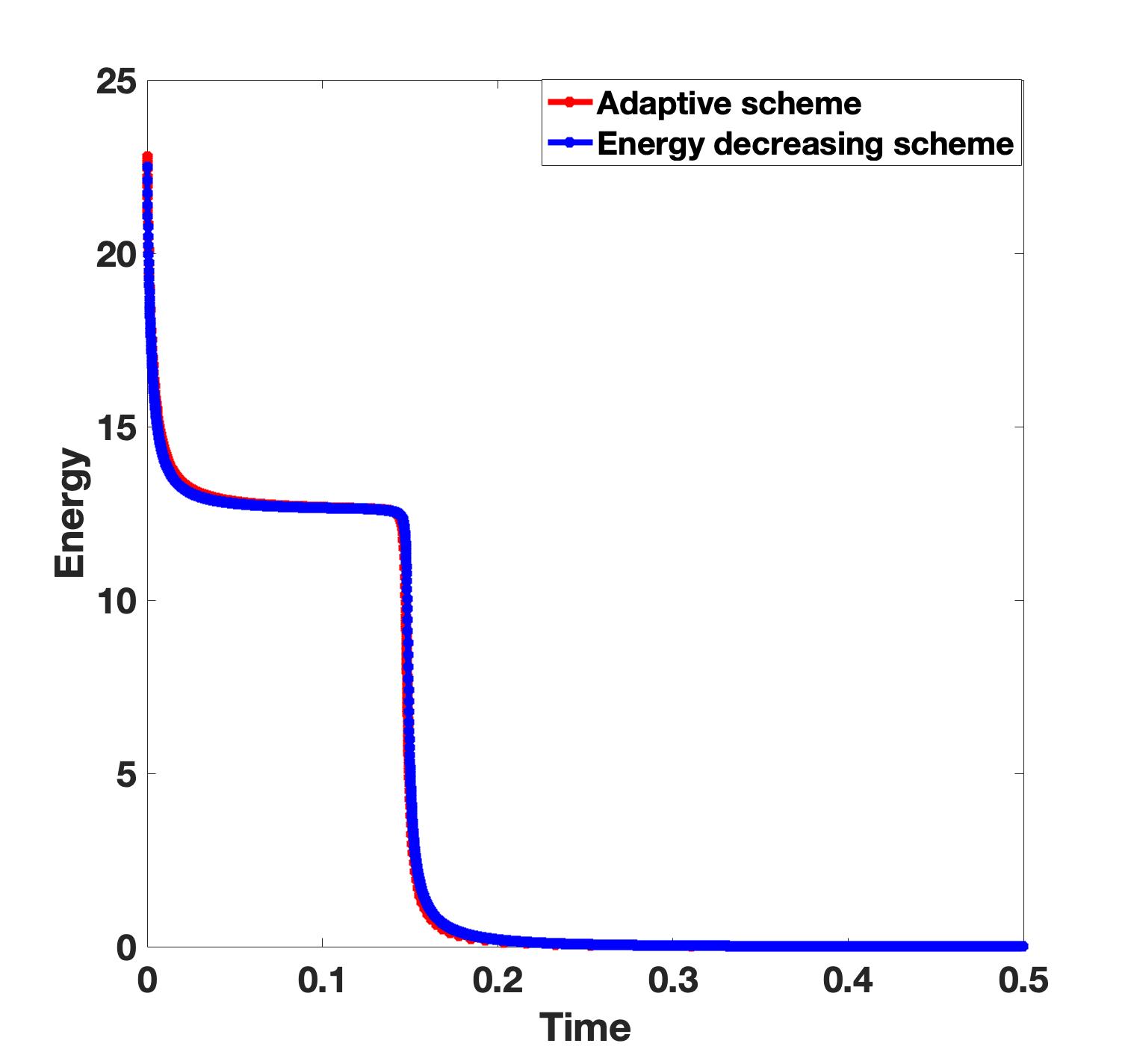}}
\subfigure[evolution of iteration numbers]{
\includegraphics[width=0.4\textwidth,clip==]{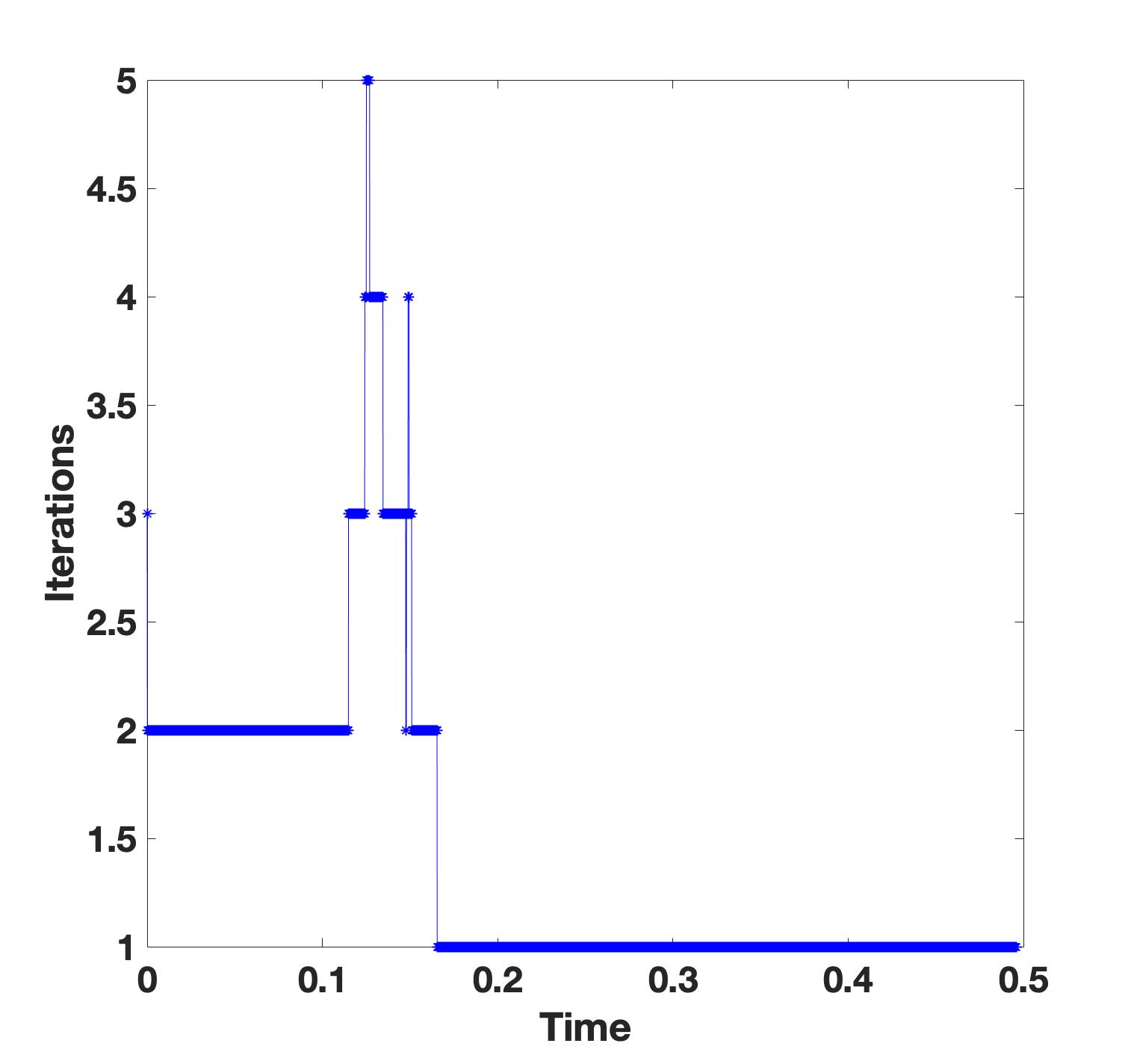}}
\caption{The evolution of  $\xi$,  energy curve and iteration numbers for Landau-Lifshitz equation \eqref{Lf:2d}  computed by \eqref{en:correction:1}-\eqref{en:correction:5} with various time step $\delta t=10^{-4}$.}\label{Fig:de:1}
\end{figure}

\begin{figure}[htbp]
\centering
\subfigure[evolution of $\xi^{n+1}$]{\includegraphics[width=0.35\textwidth,clip==]{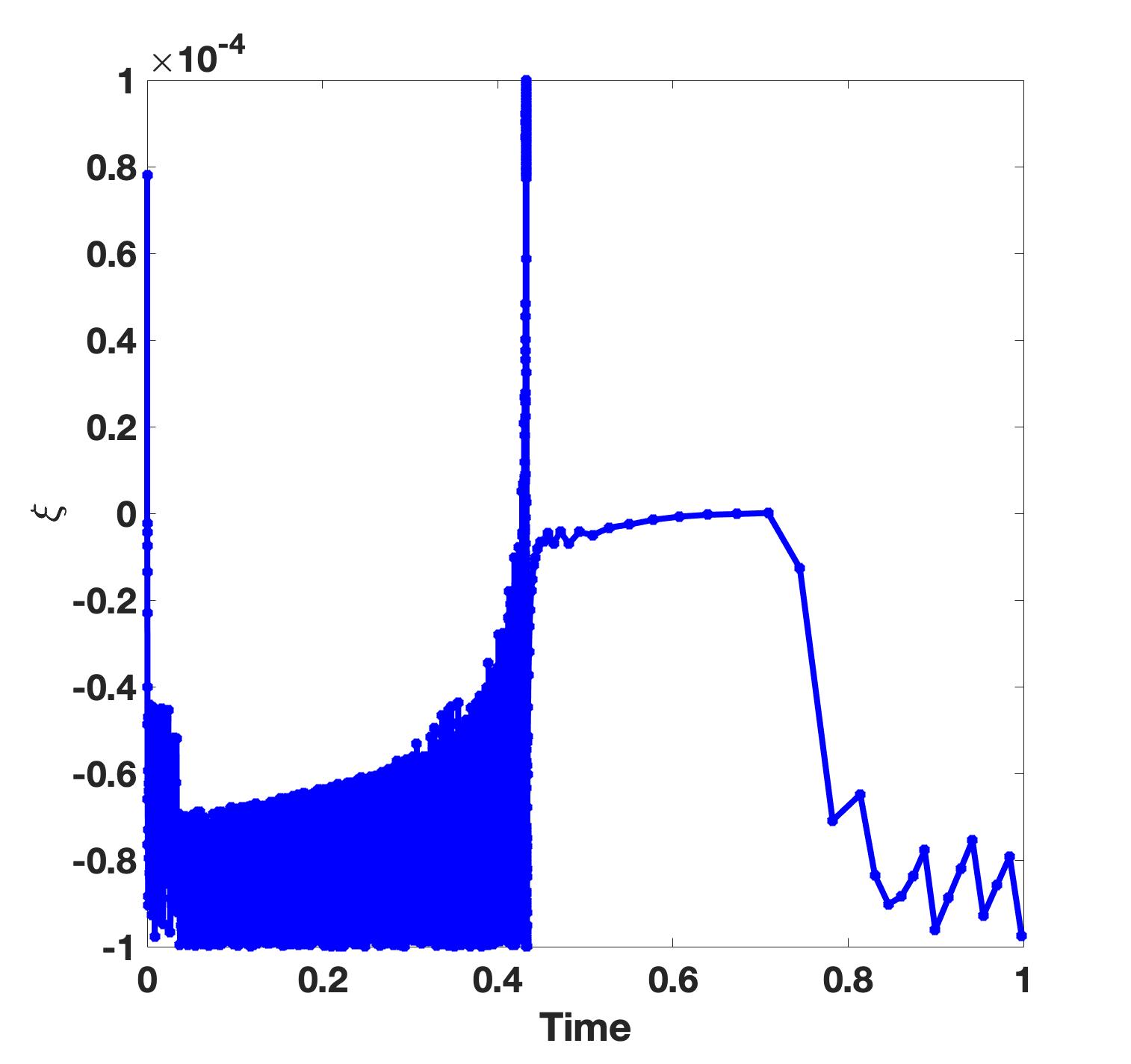}}
\subfigure[evolution of time steps]{\includegraphics[width=0.35\textwidth,clip==]{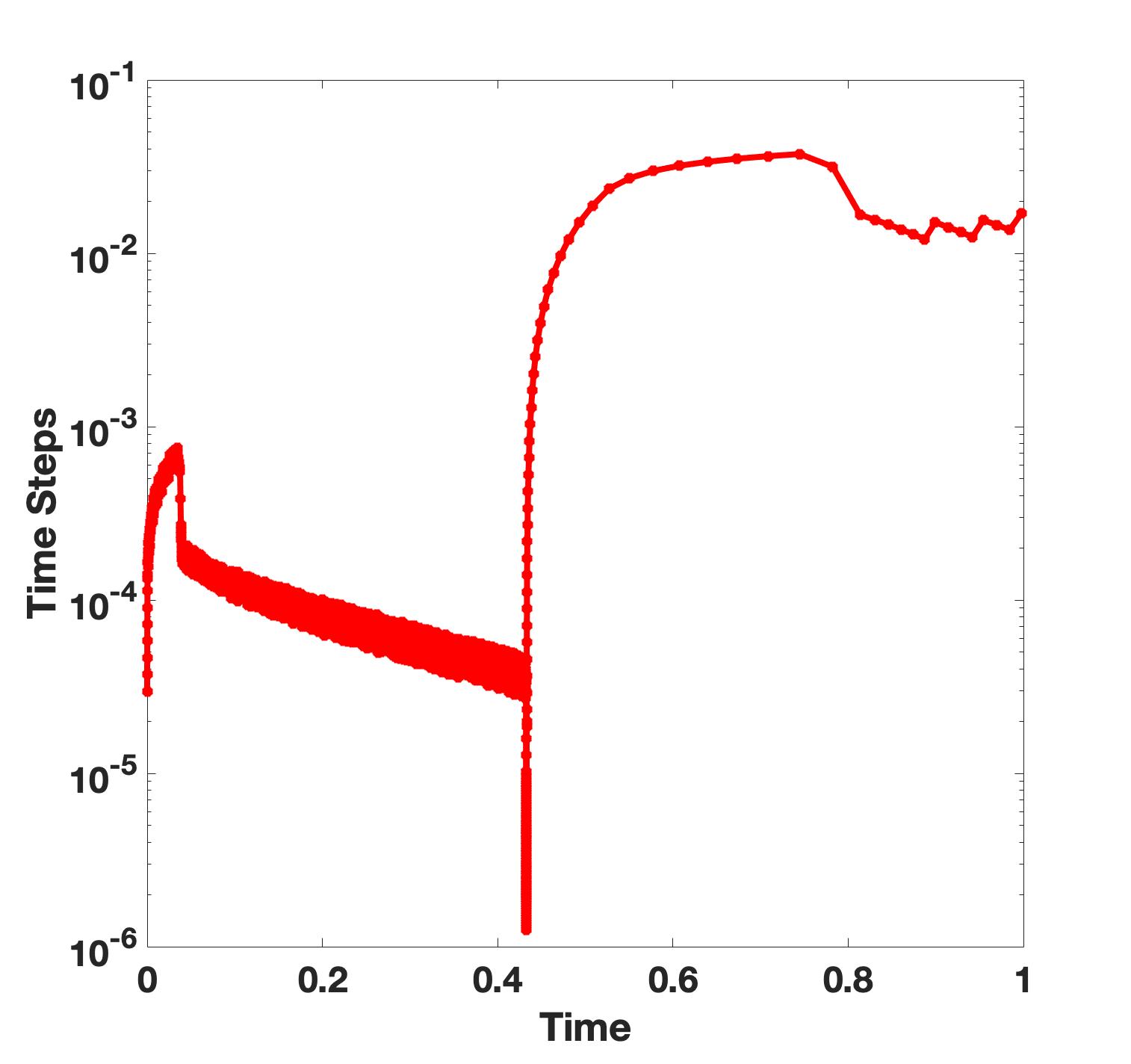}}
\subfigure[evolution of energy curves]{\includegraphics[width=0.35\textwidth,clip==]{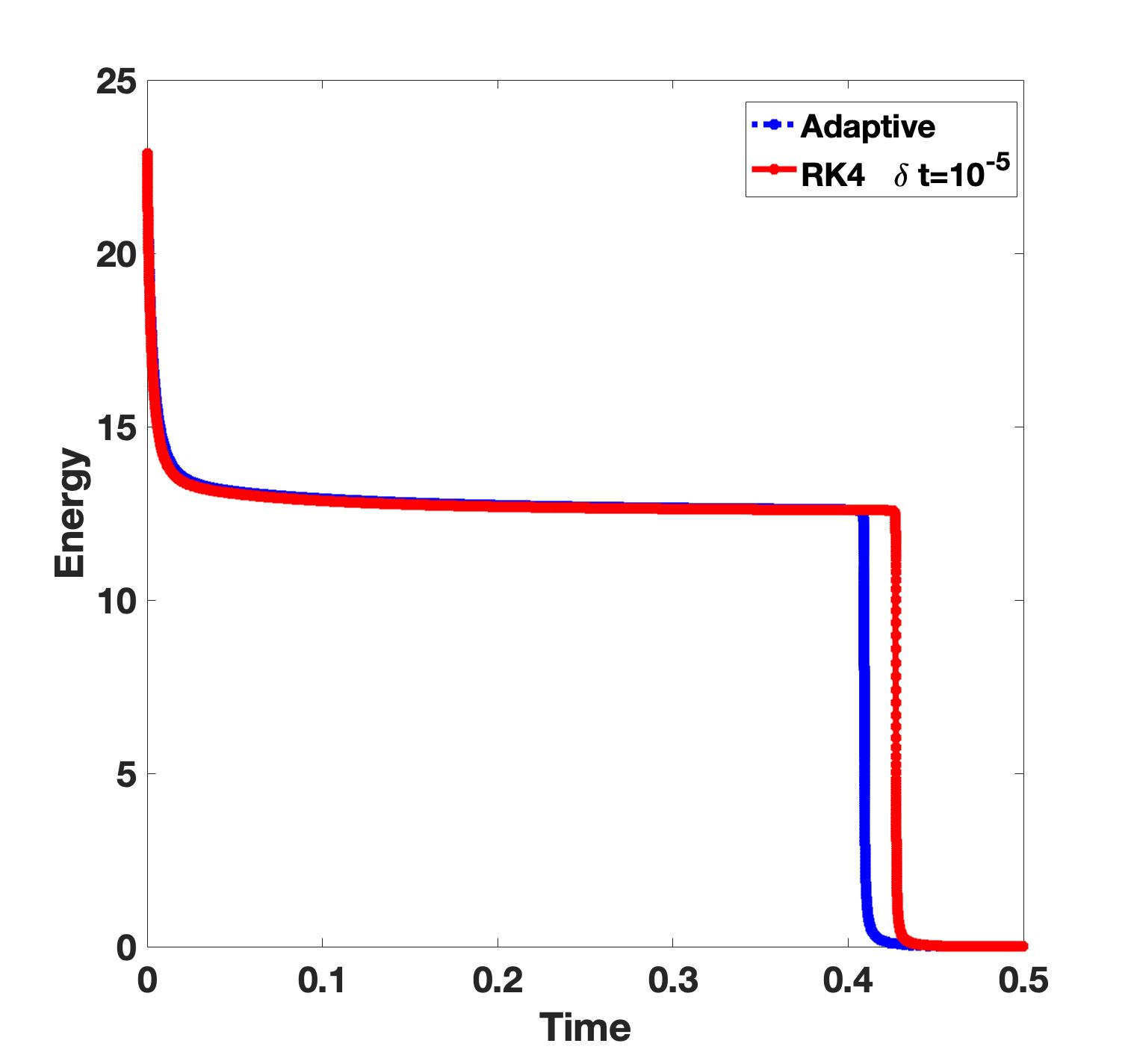}}
\subfigure[evolution of iteration numbers]{\includegraphics[width=0.35\textwidth,clip==]{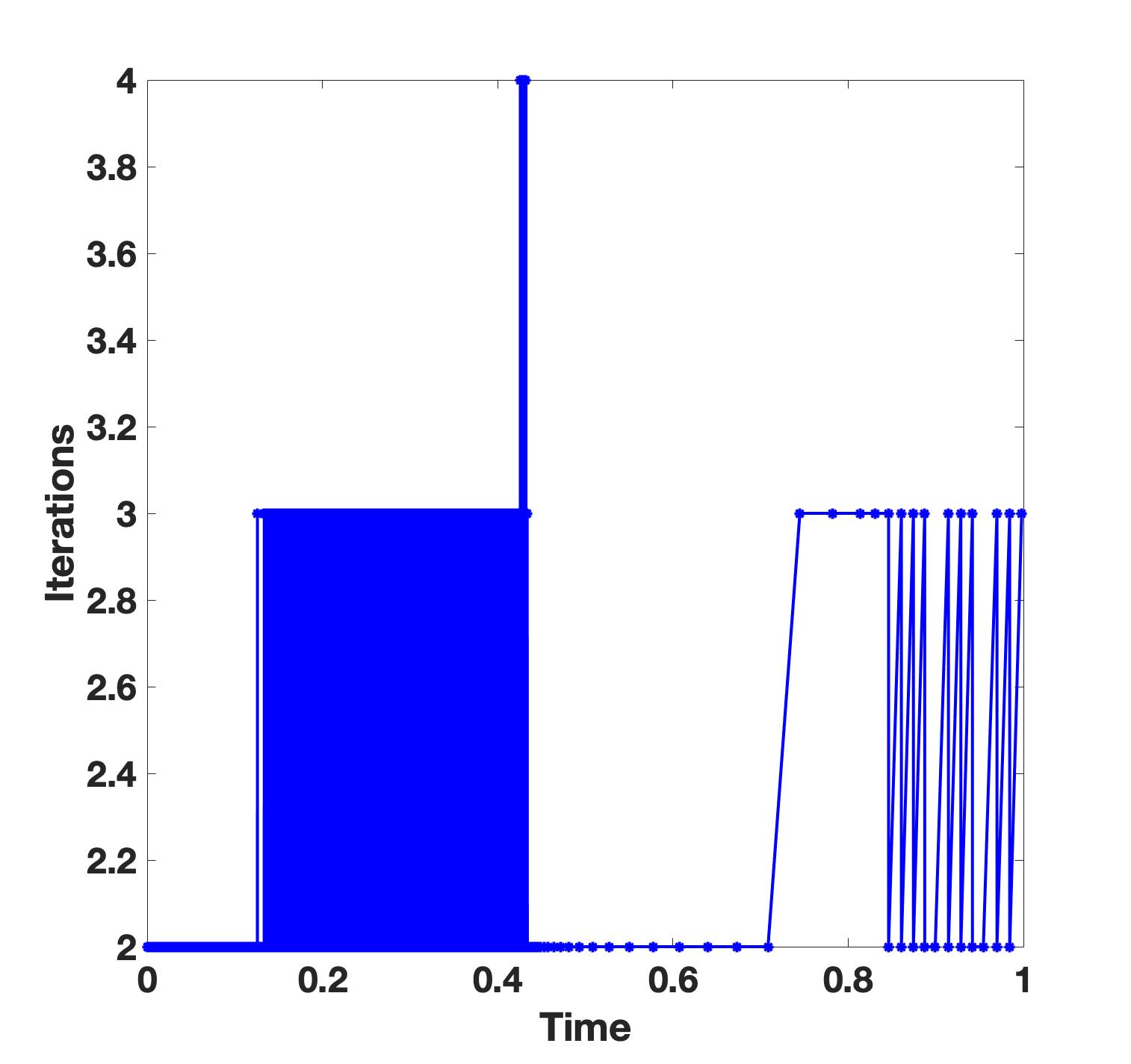}}
\caption{The evolution of  $\xi$,  energy curve and iteration numbers for Landau-Lifshitz equation  computed by \eqref{en:correction:1}-\eqref{en:correction:5} with adaptive time stepping method.}\label{adaptive:LF}
\end{figure}

First, we consider the general Landau-Lifshitz equation \eqref{Lf} with the same  initial condition  \eqref{non:ini} in the domain $\Omega=[0,2\pi)^2$ with periodic boundary conditions.   The Fourier spectral method is used in space with $128$ Fourier modes in each direction. We compute the numerical solution using energy decreasing scheme \eqref{en:correction:1}-\eqref{en:correction:5} with time step $\delta t=10^{-4}$. From Fig.\;\ref{Fig:de:1}.(a), we observe that $\xi^{n+1}$ in \eqref{en:correction:4}-\eqref{en:correction:5} is  very close to zero at most of the time and deviates from zero slightly when the energy gradient is large.   The evolution of energy computed by energy decreasing scheme \eqref{en:correction:1}-\eqref{en:correction:5}  and  its adaptive time stepping scheme  are almost the same from Fig.\;\ref{Fig:de:1}.(b).  We plot in Fig.\;\ref{Fig:de:1}.(c) the number of iteration needed to solve the nonlinear algebraic equation \eqref{en:correction:5} at each time step, and observe that  only a few iterations are needed  which indicates that the computational cost of energy decreasing scheme \eqref{en:correction:1}-\eqref{en:correction:5} is comparable with  a usual linear semi-implicit scheme.

In the second example, we test the accuracy of the adaptive time stepping scheme based on \eqref{en:correction:1}-\eqref{en:correction:5}. 
We   consider again the benchmark problem for the special Landau-Lifshitz equation \eqref{Lf:2d}  in Subsection 6.3. The results computed by scheme \eqref{en:correction:1}-\eqref{en:correction:5} with adaptive time steps are depicted in Fig.\;\ref{adaptive:LF}. We observe first that the profile of energy curve in  Fig.\;\ref{adaptive:LF}.(c) is close to  energy  curve computed by fourth-order Runge-Kutta (RK$4$) scheme  with small time step $\delta t=10^{-5}$ which indicates that we can obtain accurate numerical solutions using adaptive time stepping method.  From Fig.\;\ref{adaptive:LF}.(d), we also observe that only a few  iterations are needed at each time step to solve the nonlinear algebraic equation \eqref{en:correction:3} using the secant method. We observe from the evolution of time steps in Fig.\;\ref{adaptive:LF}.(b) that the adaptive time stepping method can significantly improve the computational efficiency.

\begin{figure}[htbp]
\centering
\includegraphics[width=0.4\textwidth,clip==]{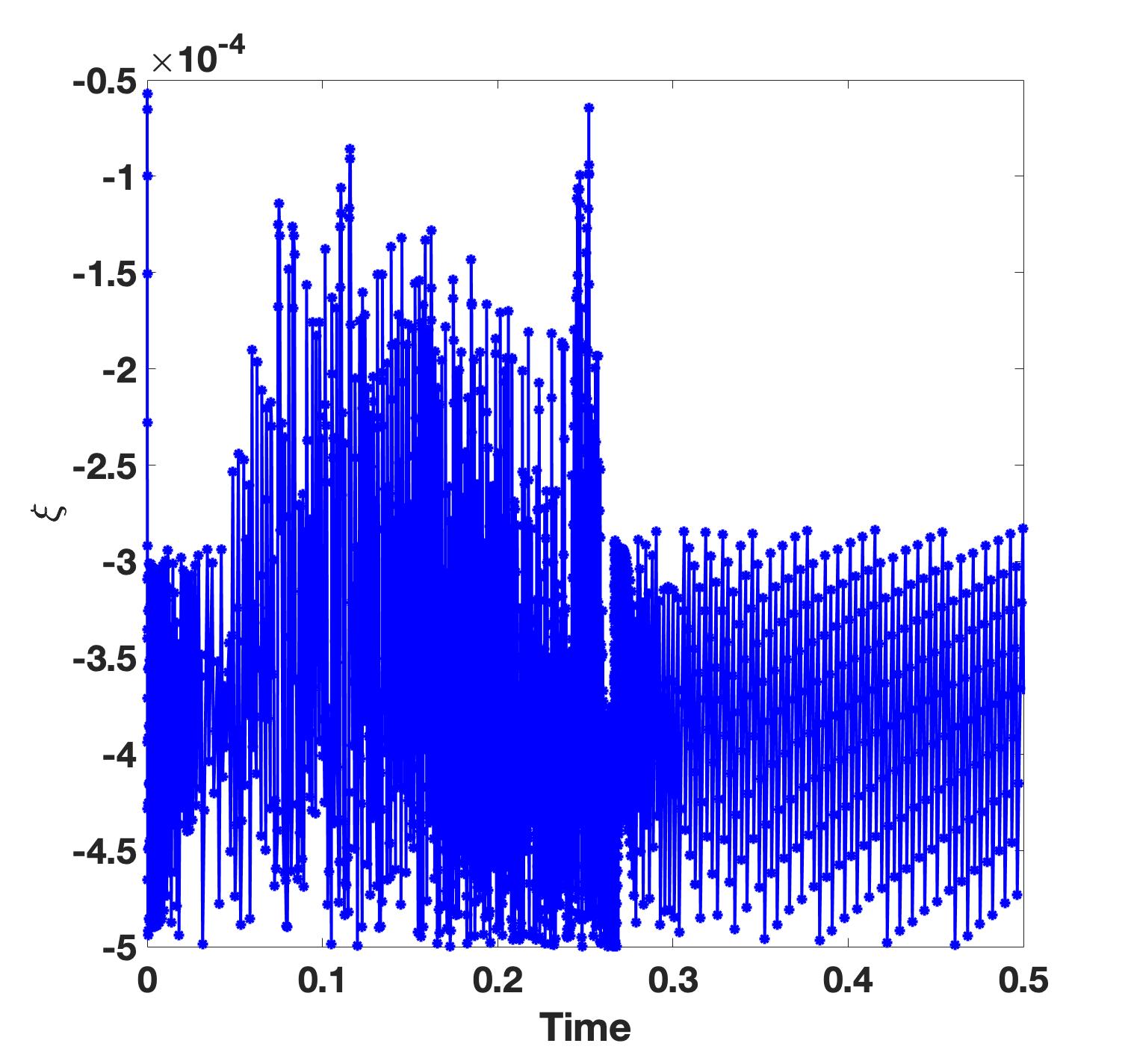}
\includegraphics[width=0.4\textwidth,clip==]{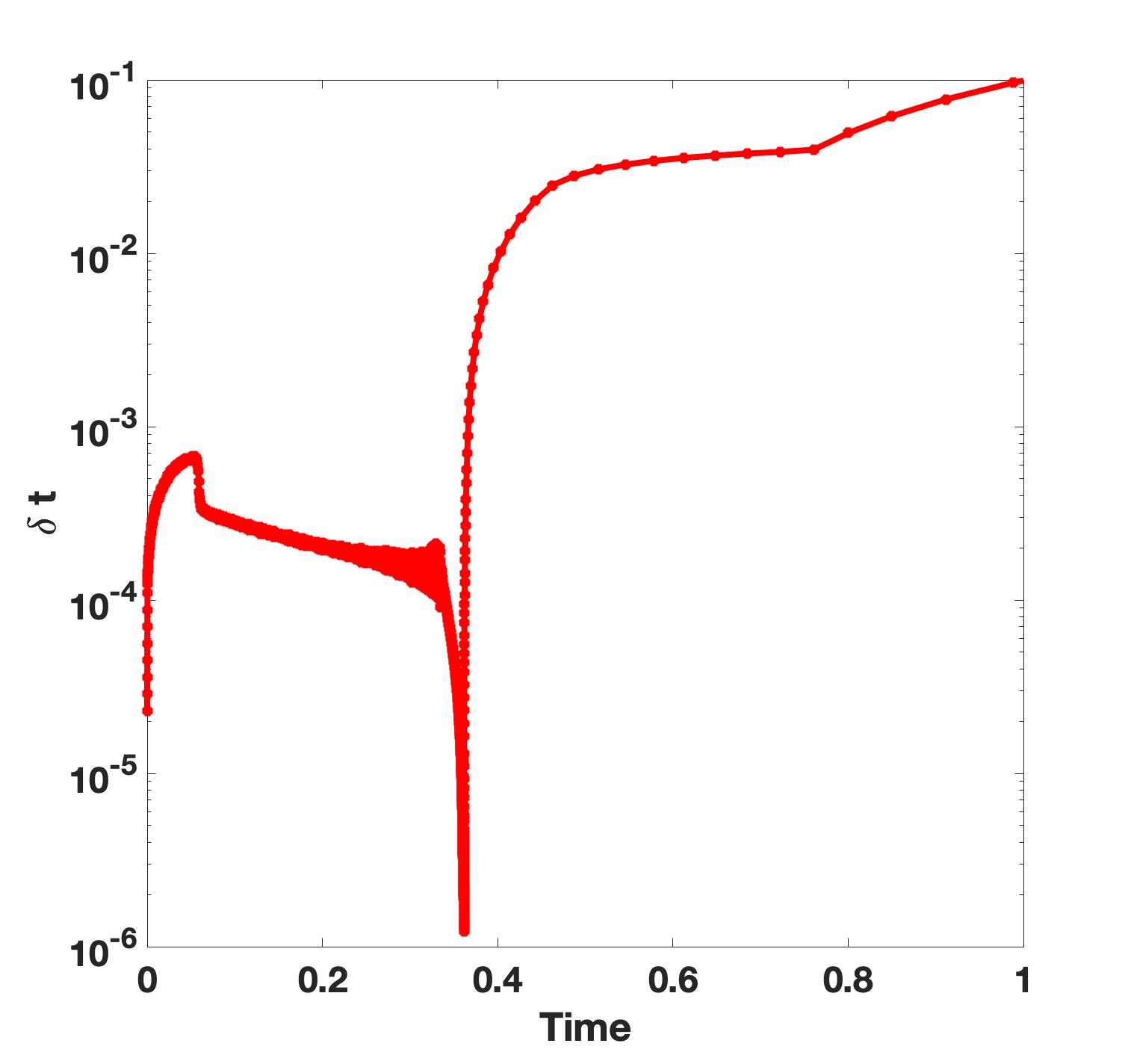}
\includegraphics[width=0.4\textwidth,clip==]{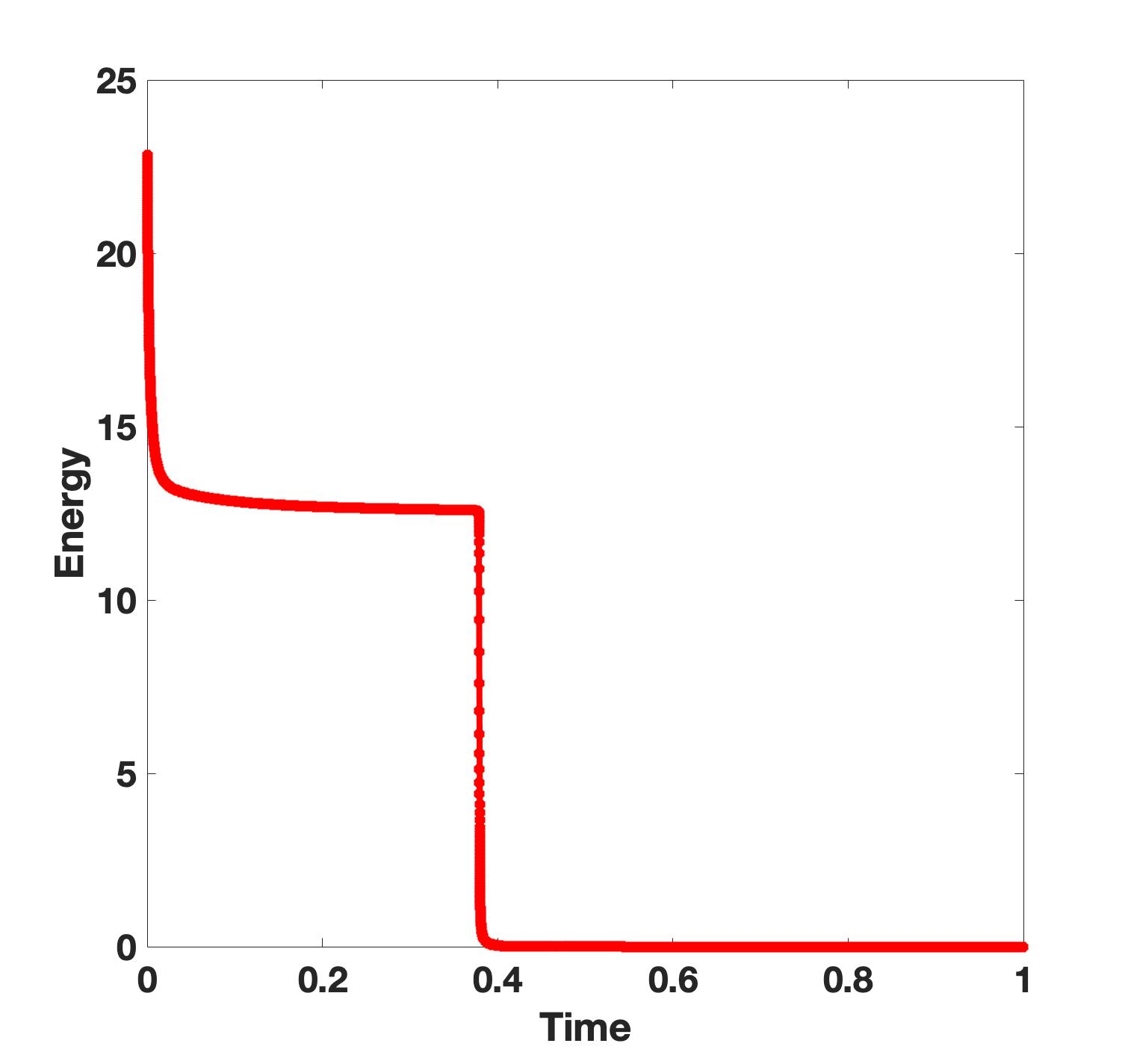}
\includegraphics[width=0.4\textwidth,clip==]{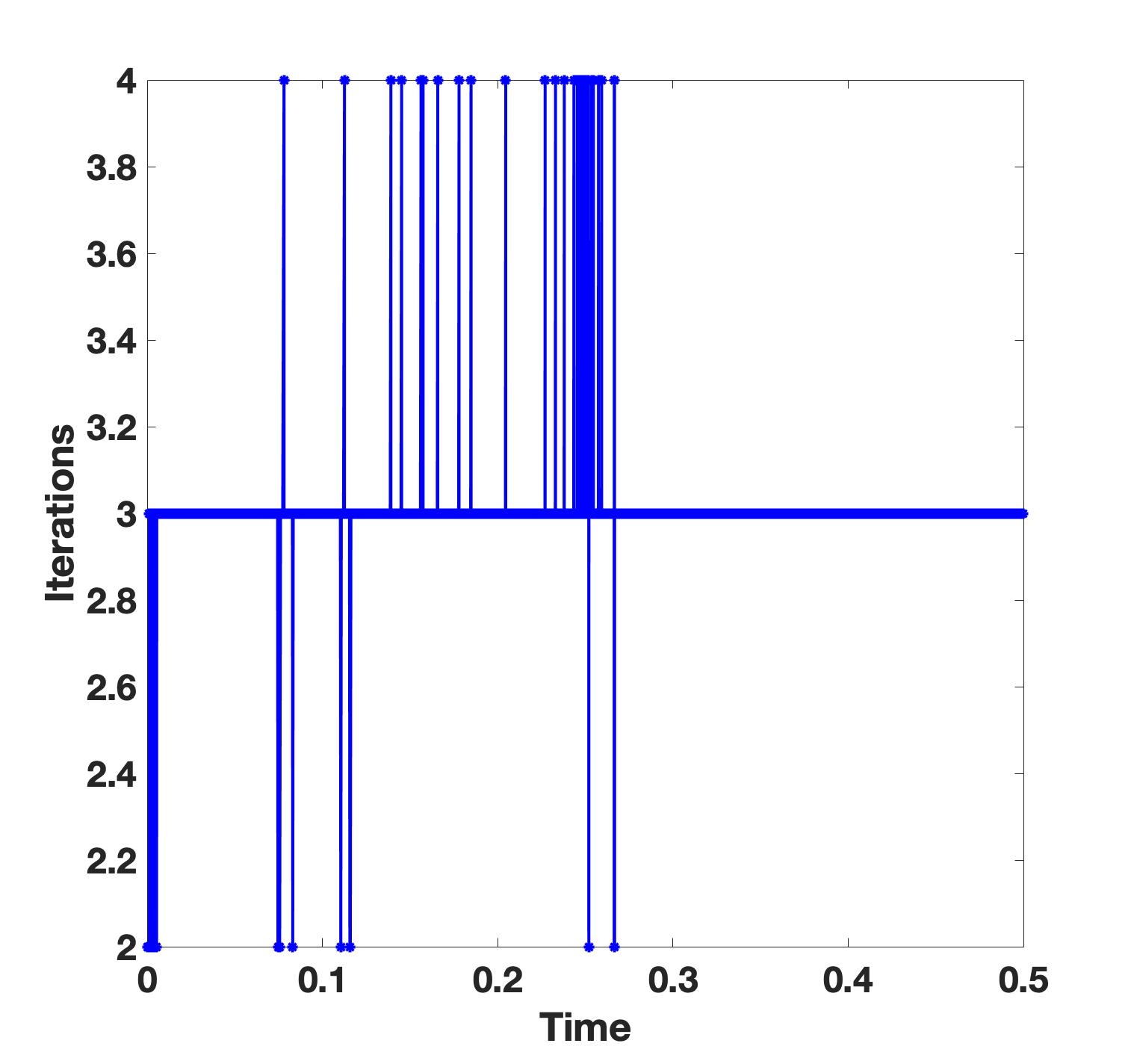}
\caption{The evolution of  $\xi$,  energy curve and iteration numbers for the general Landau-Lifshitz equation  computed by \eqref{seidel:LLG:correction:2}-\eqref{seidel:LLG:correction:5} with adaptive time stepping method.}\label{adaptive:LLG}
\end{figure}

In the third example, we  solve the Laudau-Lifshitz equation \eqref{Lf} with $\beta=\gamma=1$ in the domain $[-\frac 12,\frac 12)^2$  using the adaptive time stepping scheme with tolerance $ |\xi|\leq tol=5\times 10^{-5}$  based on \eqref{seidel:LLG:correction:2}-\eqref{seidel:LLG:correction:5} and $64^2$ Fourier modes in space. The initial condition is also chosen to be \eqref{non:ini}.  We observe from Fig.\;\ref{adaptive:LLG} that larger time steps can be used when the energy curve change slowly, and the computed energy  decreases with time.

\begin{comment}
\begin{figure}[htbp]
\centering
\includegraphics[width=0.4\textwidth,clip==]{adaptive_xi.jpg}
\includegraphics[width=0.4\textwidth,clip==]{adaptive_steps.jpg}
\includegraphics[width=0.4\textwidth,clip==]{adaptive_en2.jpg}
\includegraphics[width=0.4\textwidth,clip==]{adaptive_iters.jpg}
\caption{The evolution of  $\xi$,  energy curve and iteration numbers for Landau-Lifshitz equation  computed by \eqref{en:correction:1}-\eqref{en:correction:5} with adaptive time stepping method.}\label{Adaptive:Fig:de:1}
\end{figure}
\end{comment}

\subsection{phenomenon of blowup}

\begin{figure}[htbp]
\centering
\subfigure[$m$ at $t=0$.]{
\includegraphics[width=0.32\textwidth,clip==]{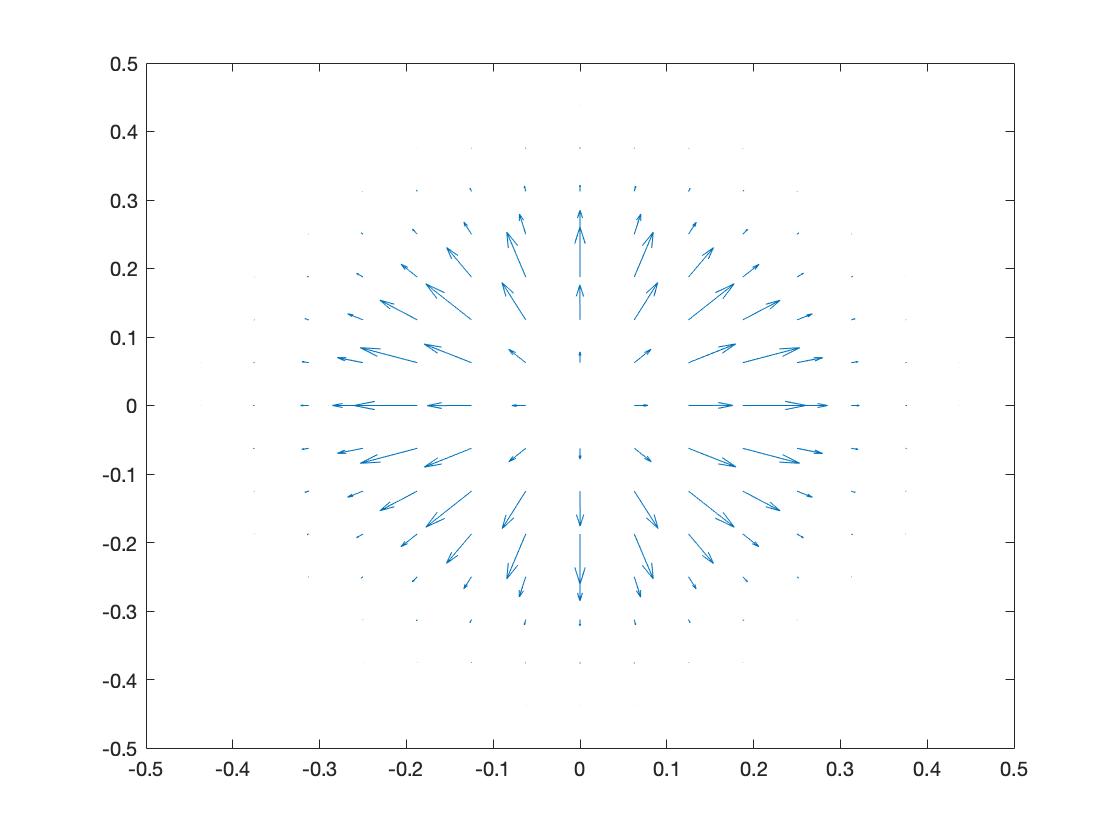}}
\subfigure[$m$ at $t=0.001$.]{
\includegraphics[width=0.32\textwidth,clip==]{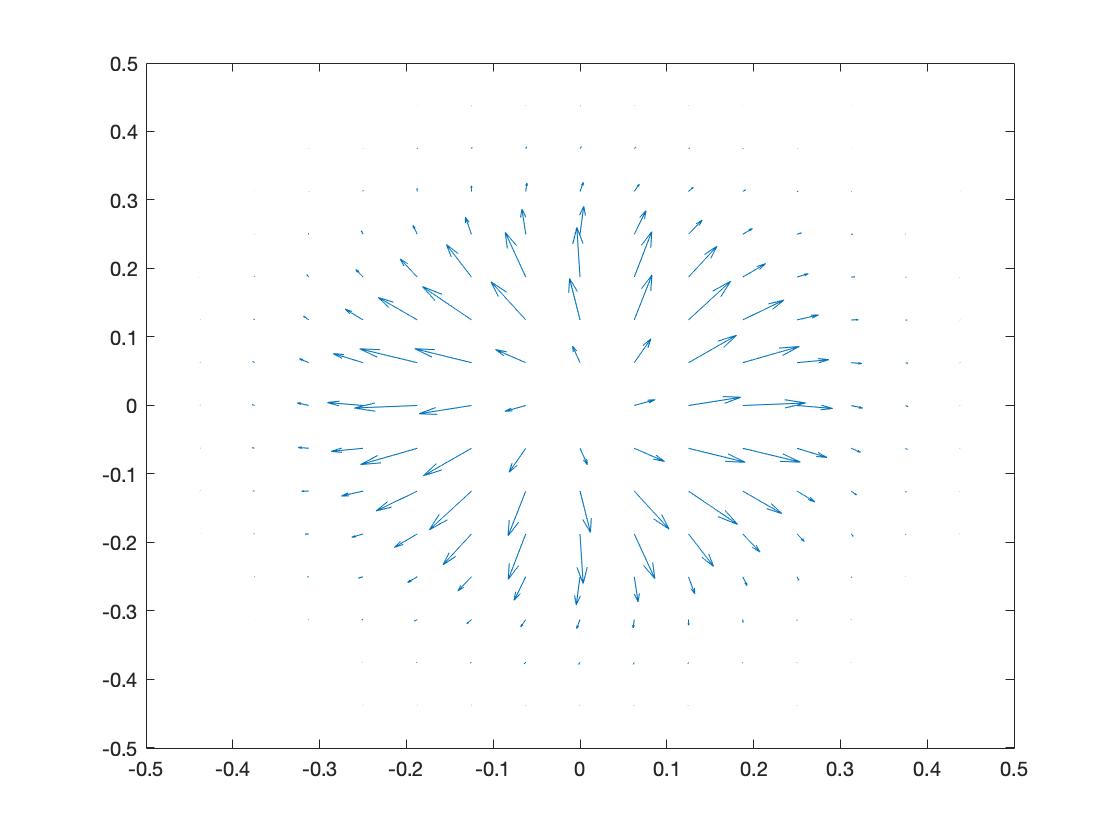}}
\subfigure[$m$ at $t=0.01$.]{
\includegraphics[width=0.32\textwidth,clip==]{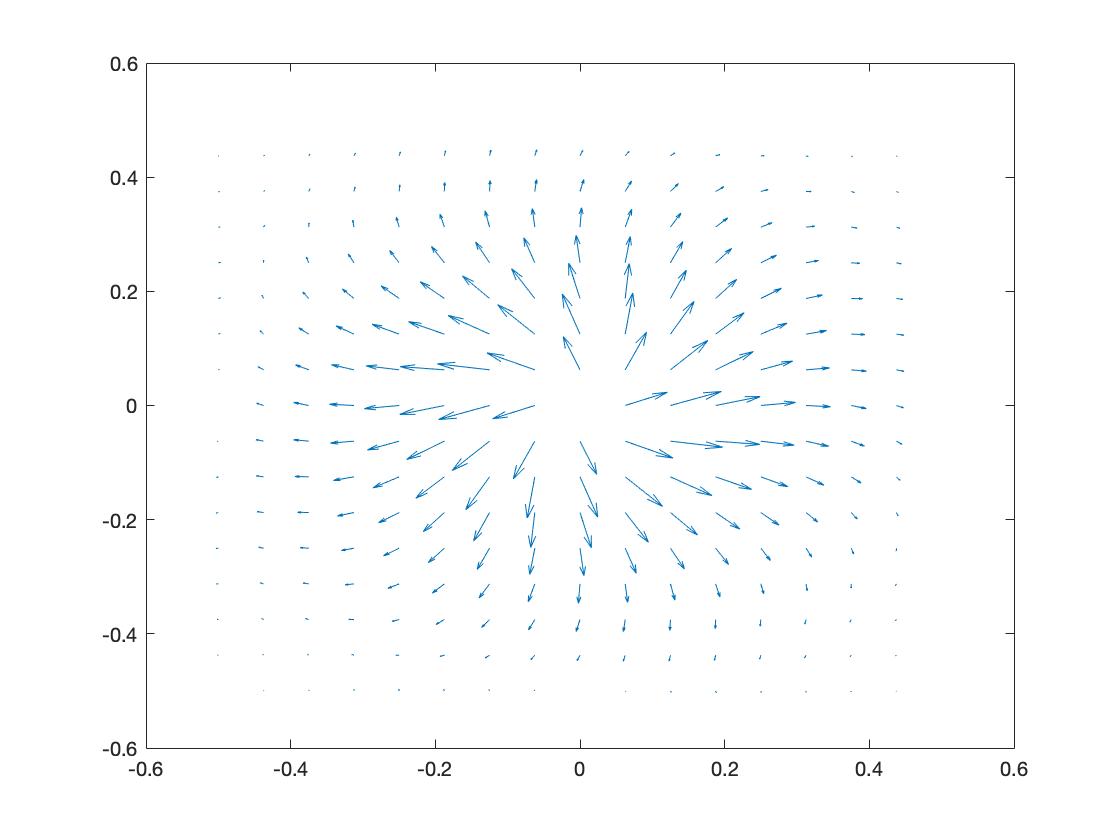}}
\subfigure[$m$ at $t=0.05$.]{
\includegraphics[width=0.32\textwidth,clip==]{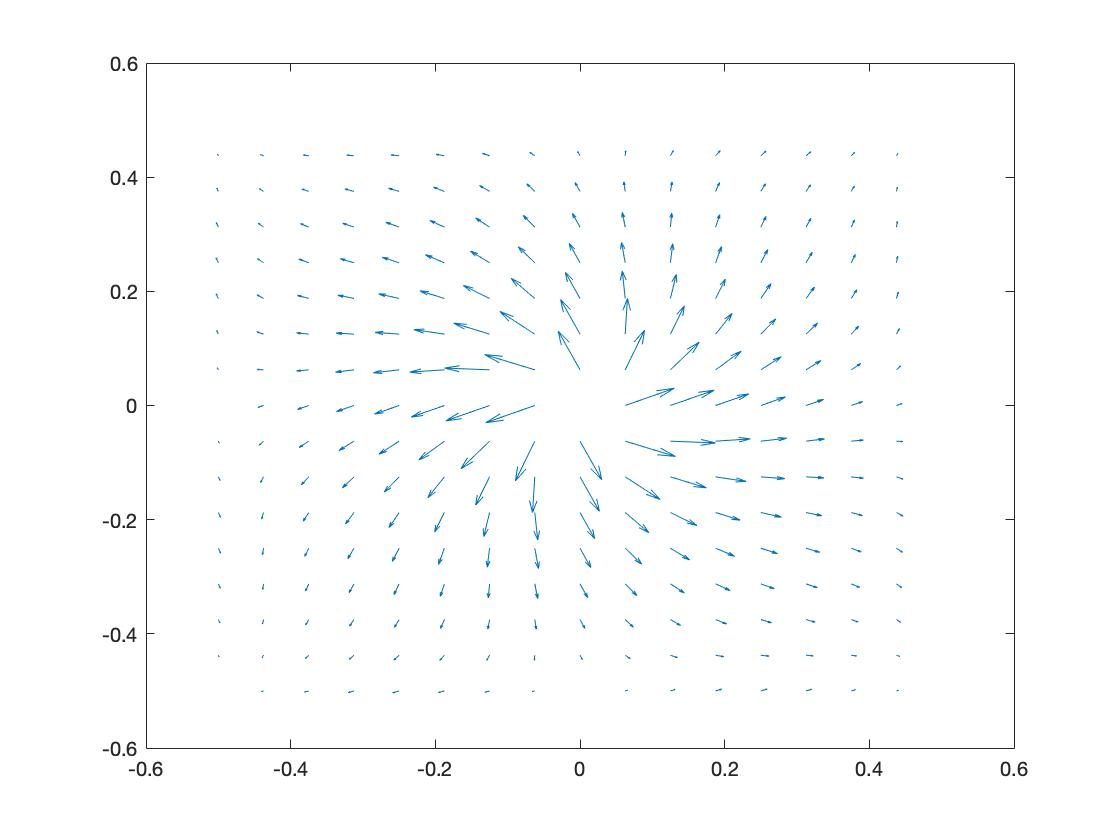}}
\subfigure[$m$ at $t=0.1$.]{
\includegraphics[width=0.32\textwidth,clip==]{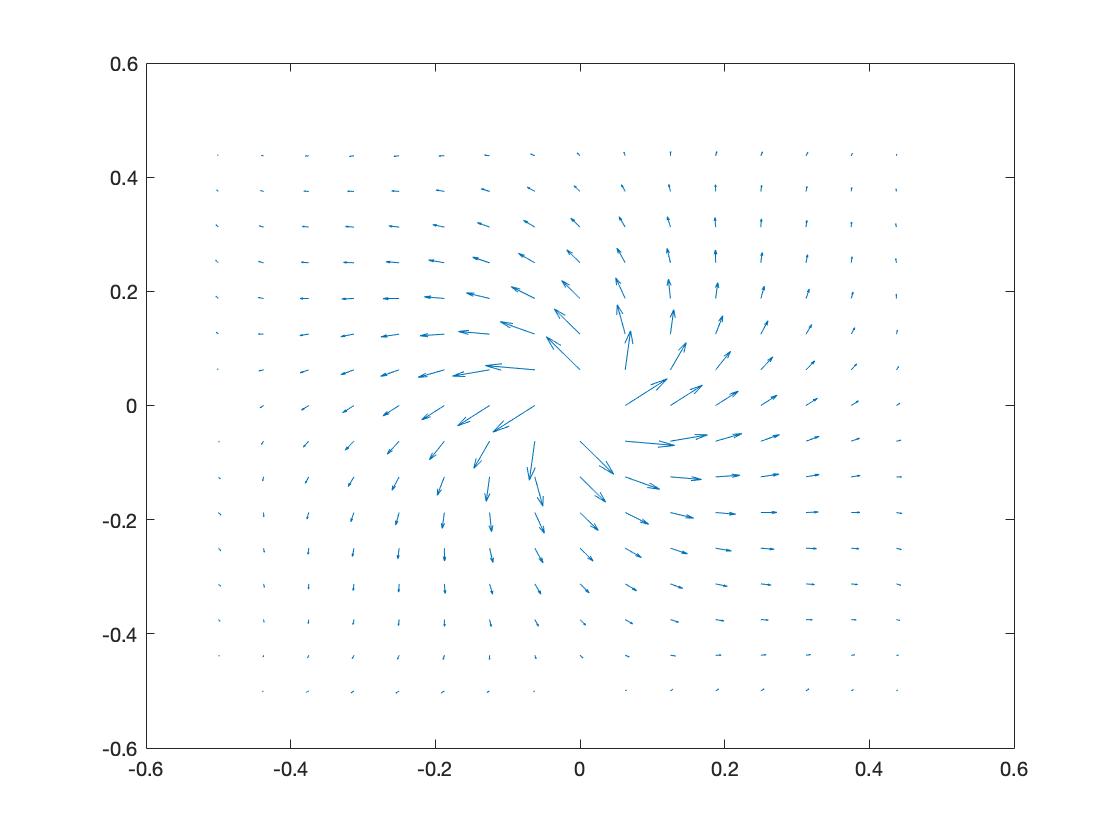}}
\subfigure[$m$ at $t=0.2$]{
\includegraphics[width=0.32\textwidth,clip==]{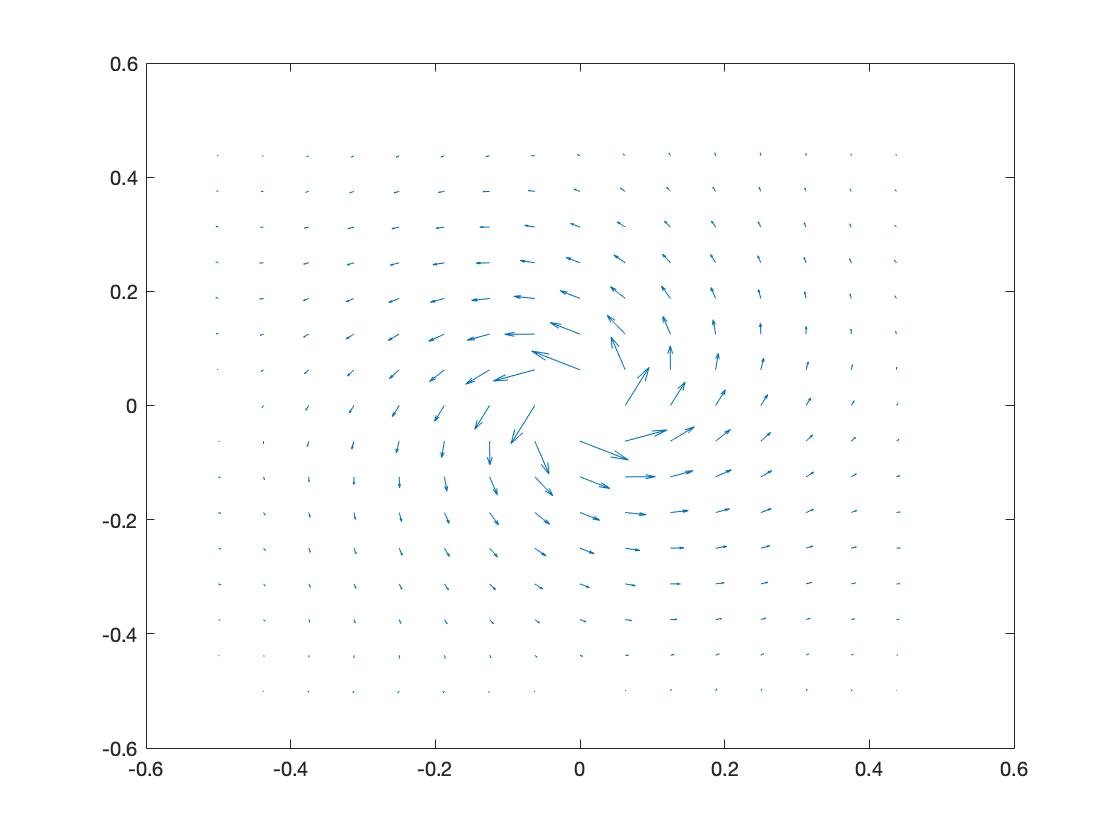}}
\subfigure[$m$ at $t=0.4$.]{
\includegraphics[width=0.32\textwidth,clip==]{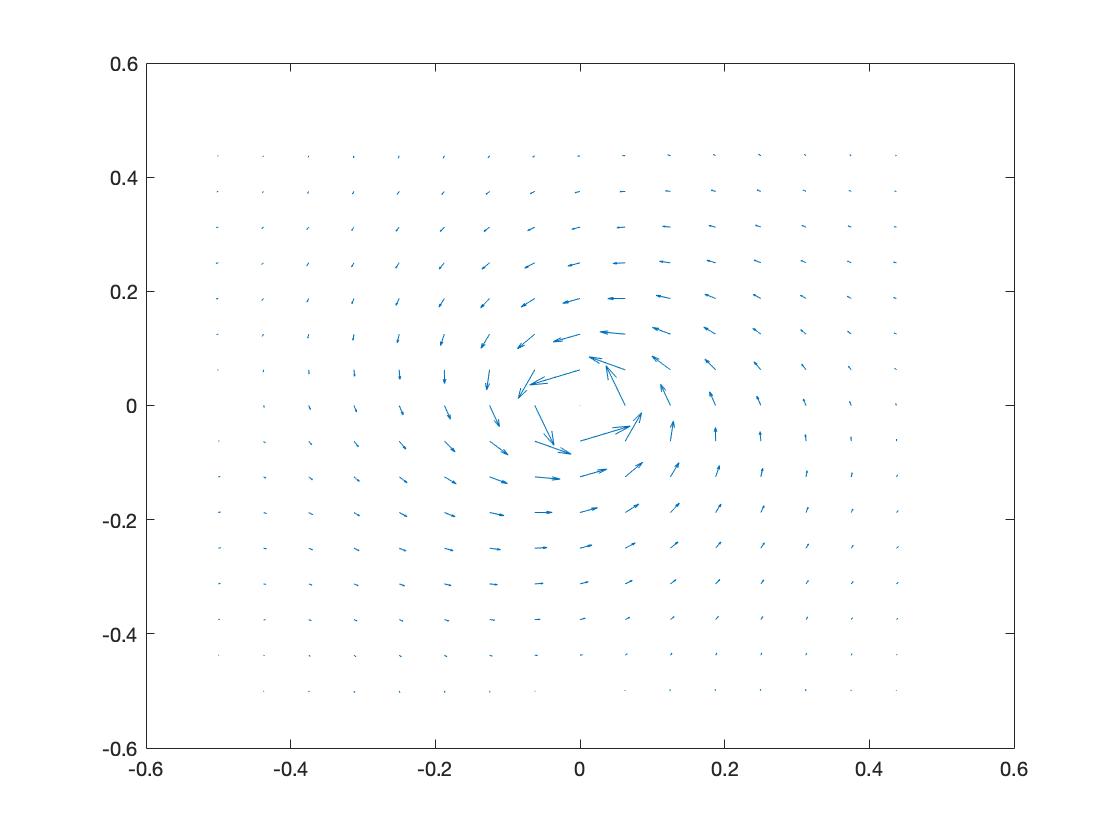}}
\subfigure[$m$ at $t=0.5$.]{
\includegraphics[width=0.32\textwidth,clip==]{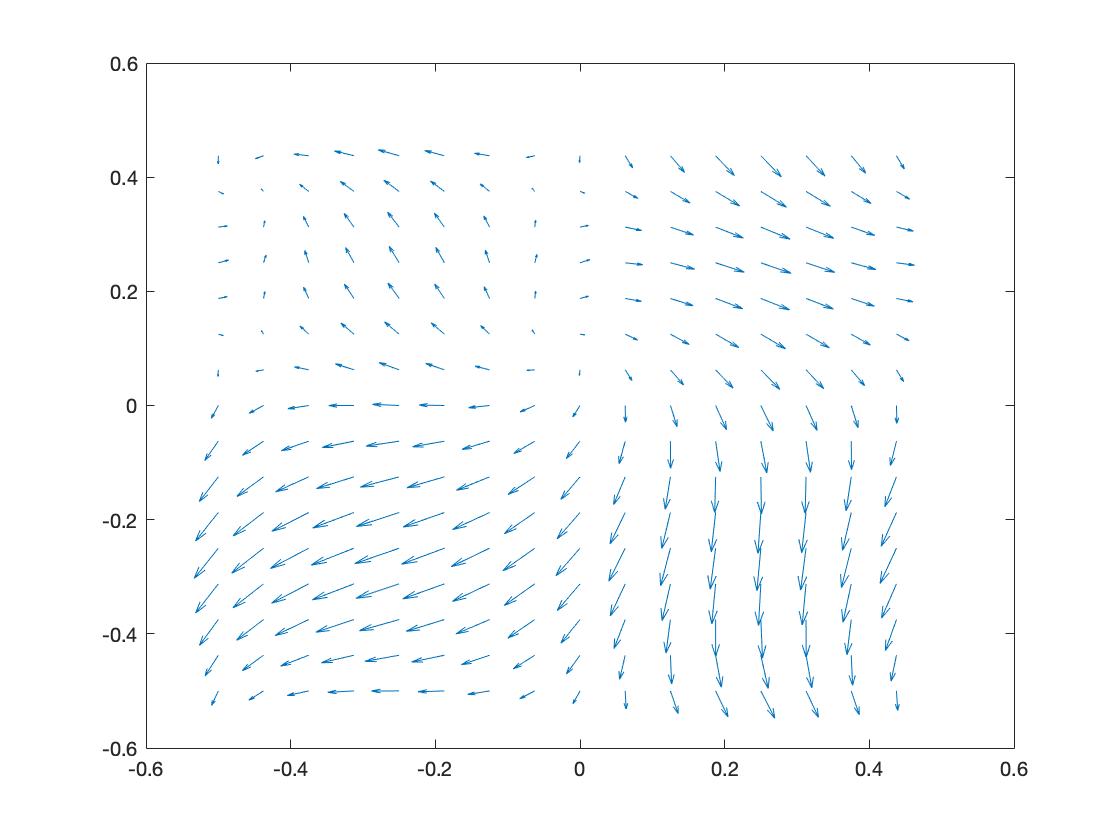}}
\subfigure[$m$ at $t=0.6$.]{
\includegraphics[width=0.32\textwidth,clip==]{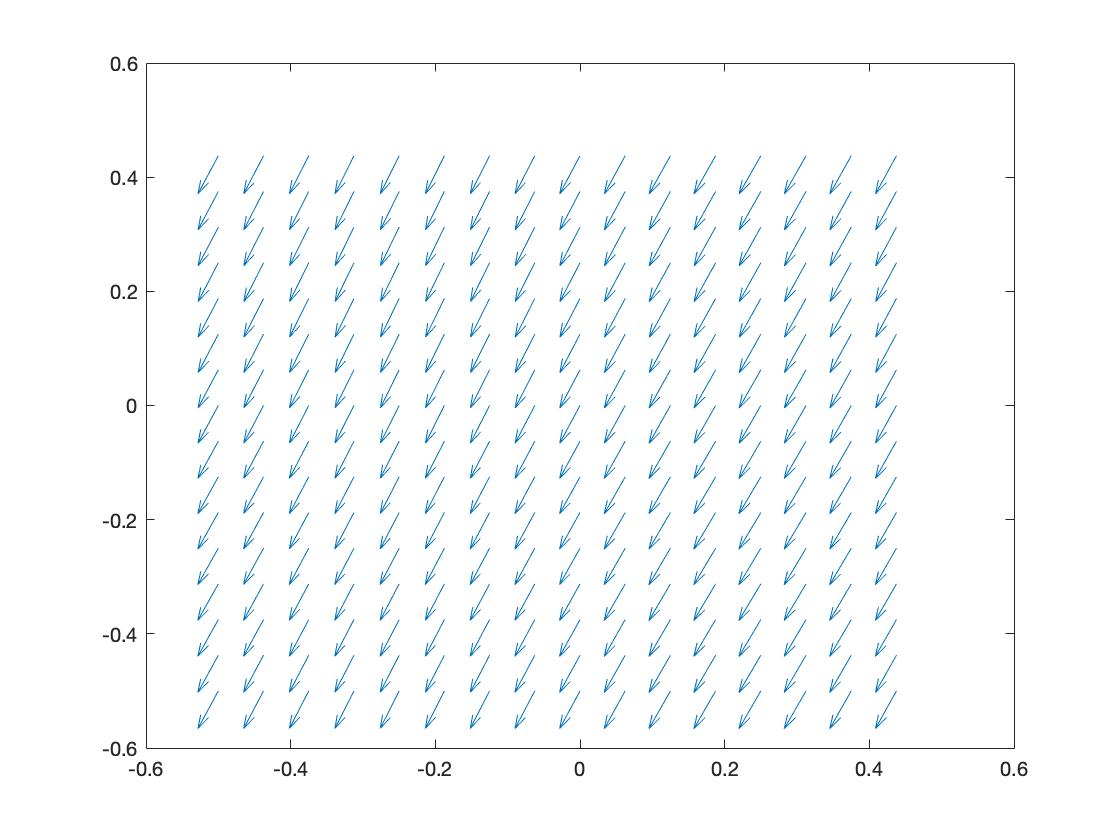}}
\caption{(a)-(i): Numerical solutions $m$ at $t=0, 0.001, 0.01, 0.05, 0.1, 0.2, 0.4, 0.5, 0.6$ projected on $xy$-plane using the second-order Gauss-Seidel  predictor-correction scheme \eqref{LLG:correction:1b}-\eqref{LLG:correction:3b}  with $\delta t=10^{-5}$.}\label{Fourier}
\end{figure}

\begin{figure}[htbp]
\centering
\subfigure[$m$ at $t=0$.]{
\includegraphics[width=0.32\textwidth,clip==]{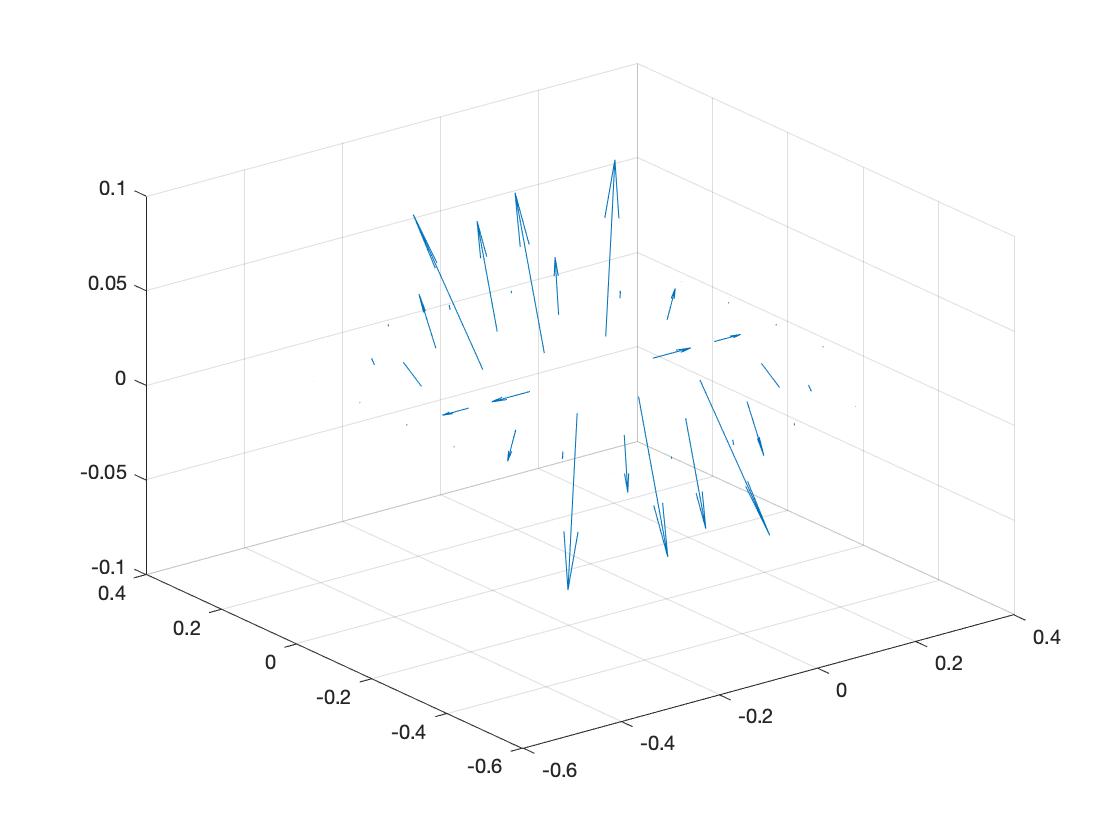}}
\subfigure[$m$ at $t=0.4$.]{
\includegraphics[width=0.32\textwidth,clip==]{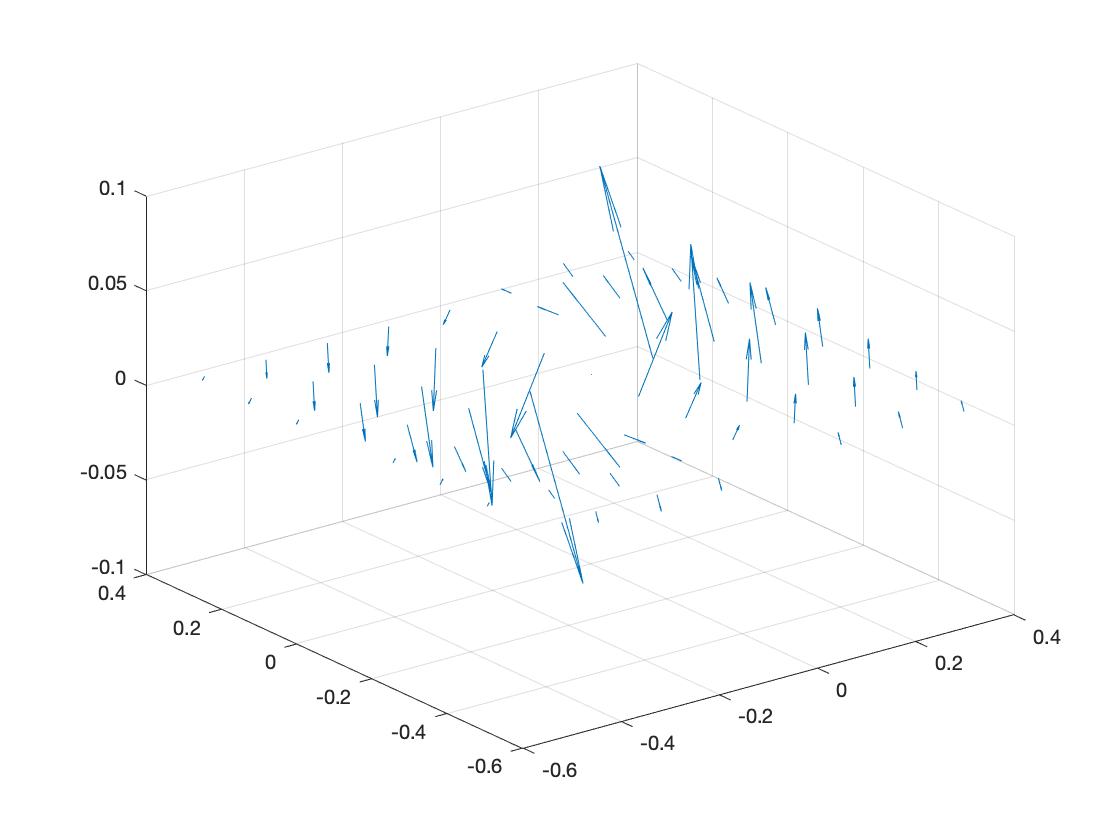}}
\subfigure[$m$ at $t=0.6$.]{
\includegraphics[width=0.32\textwidth,clip==]{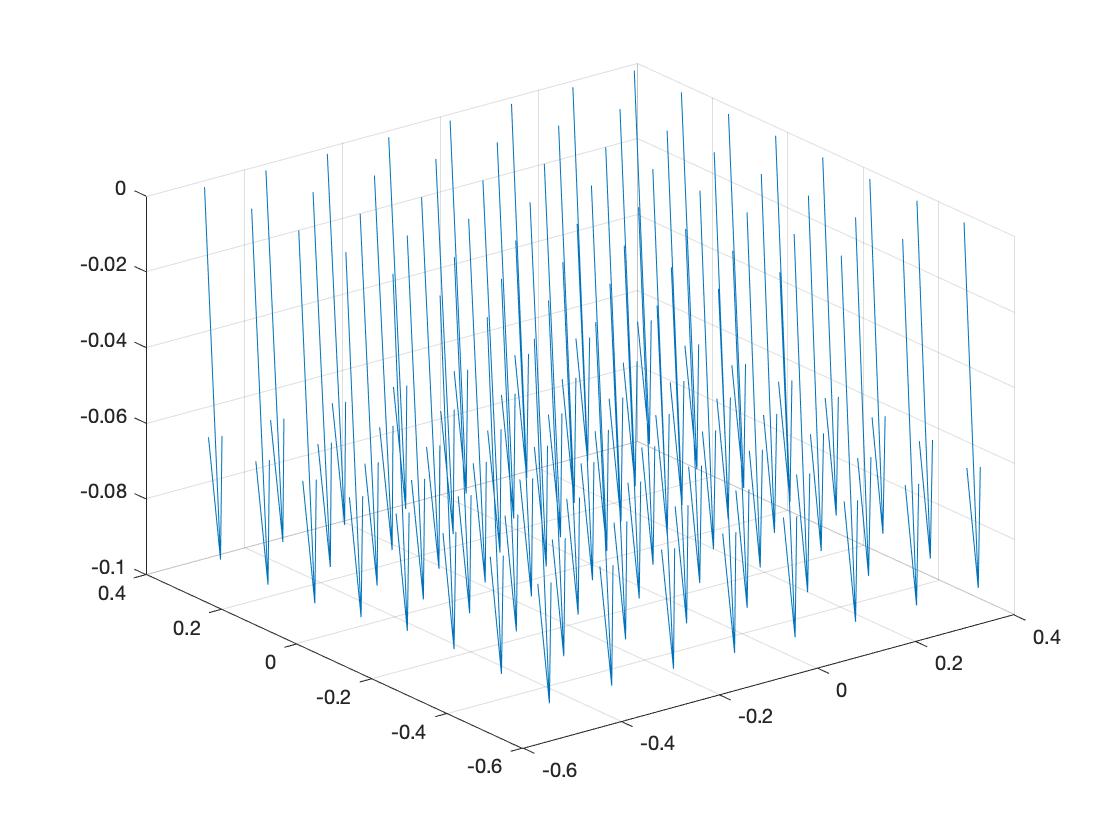}}
\caption{(a)-(c): Numerical solutions $m$ at $t=0, 0.4, 0.6$ using the second-order Gauss-Seidel  predictor-correction scheme \eqref{LLG:correction:1b}-\eqref{LLG:correction:3b} with $\delta t=10^{-5}$.}\label{Fourier_3d}
\end{figure}

\begin{figure}[htbp]
\centering
\subfigure[evolution of energy curves]{
\includegraphics[width=0.4\textwidth,clip==]{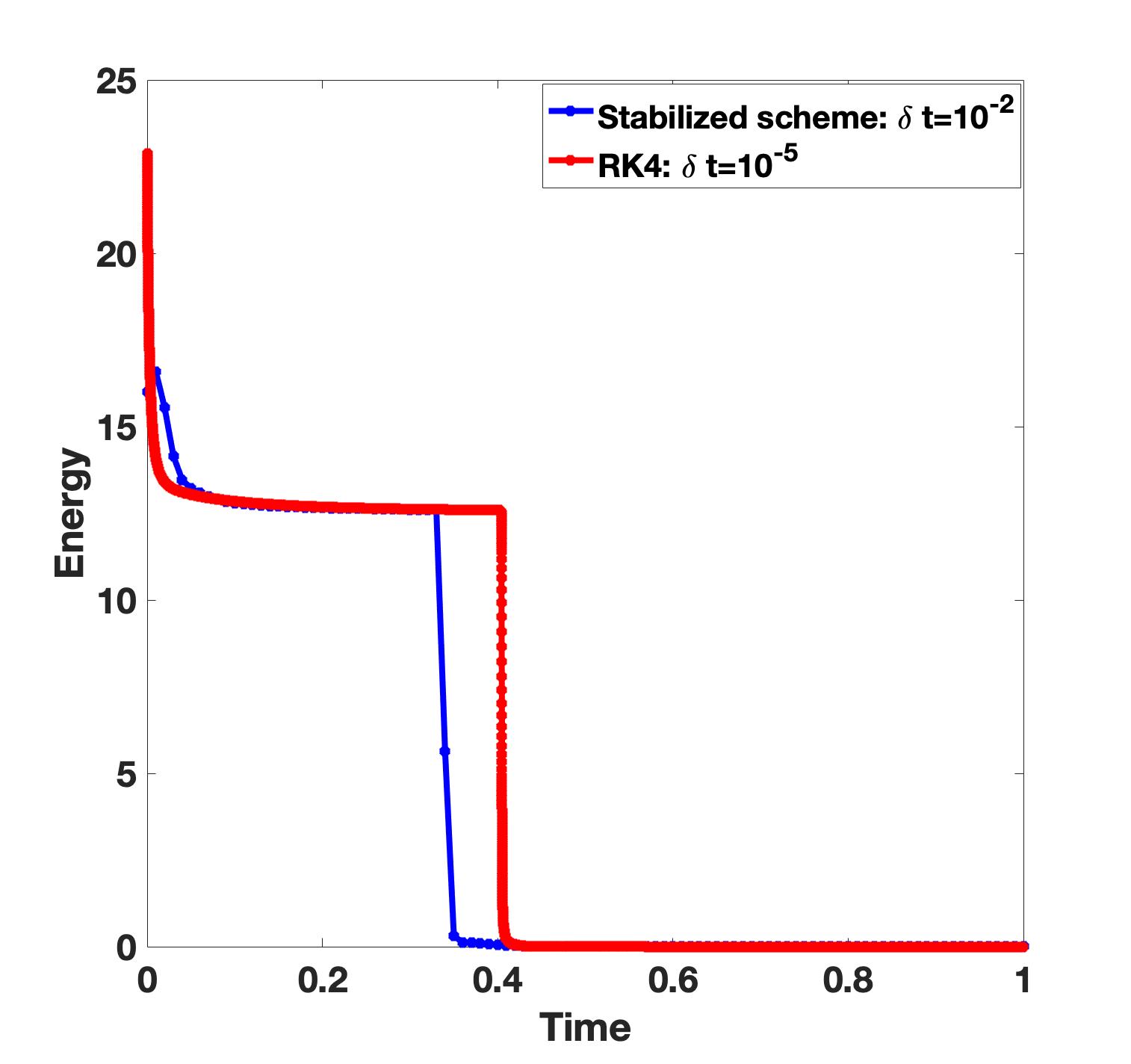}}
\subfigure[evolution of $|m|_{W^{1,\infty}}$ ]{
\includegraphics[width=0.4\textwidth,clip==]{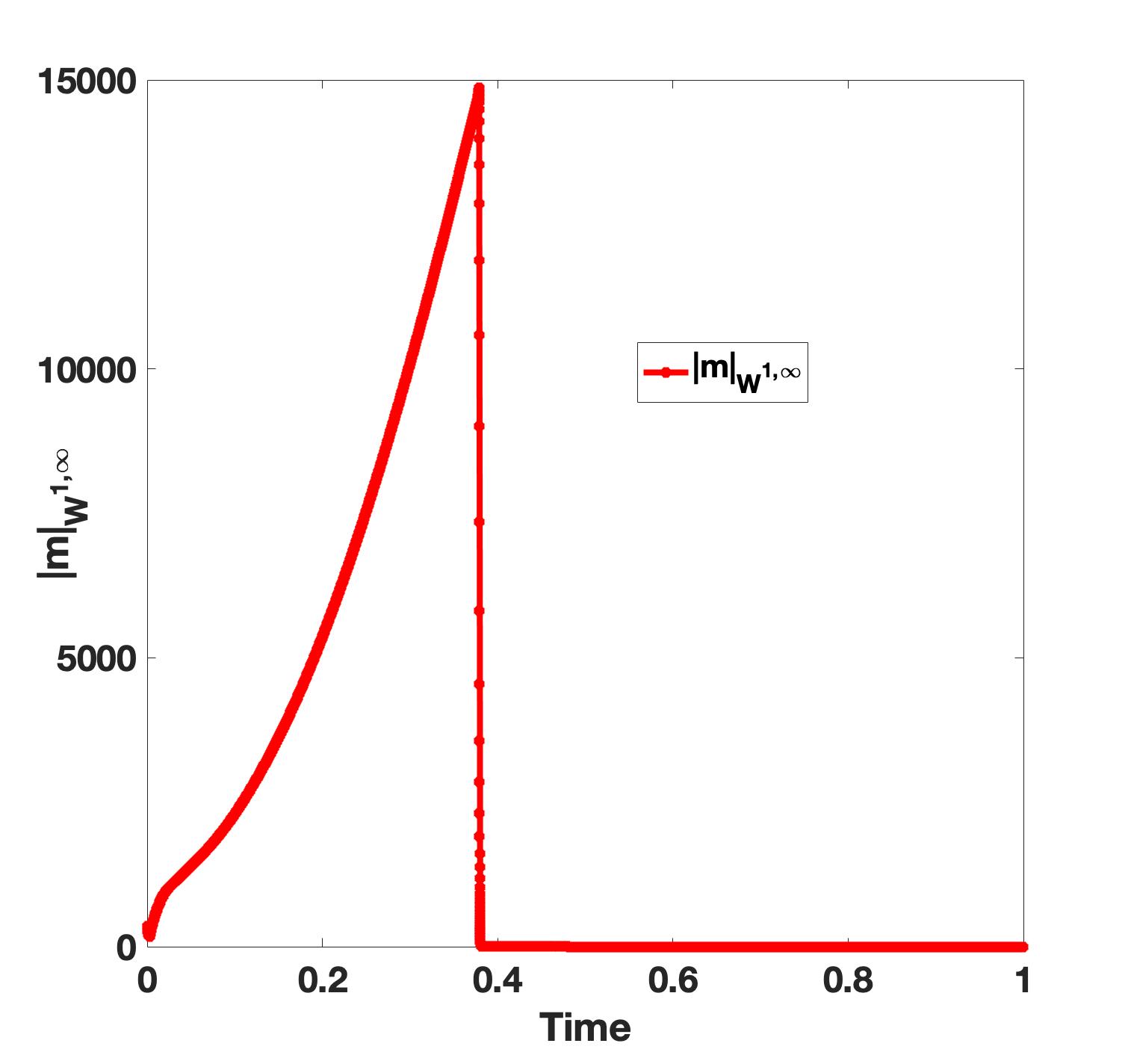}}
\caption{(a). Comparison of energy curves using stabilized type-II Gauss-Seidel scheme  with stabilized constant $S=0.5$. (b). The evolution of  $\|m\|_{W^{1,\infty}}$ for the general Landau-Lifshitz equation  computed by Gauss-Seidel \eqref{LLG:correction:1b}-\eqref{LLG:correction:3b} with time step $\delta t=10^{-5}$.}\label{grad_m}
\end{figure}

In this subsection,  we investigate the possible blowup of the general Landau-Lifshitz equation \eqref{Lf} with certain smooth initial condition  which has been studied in \cite{bartels2008numerical,an2021optimal}. We set $\Omega =[-\frac 12,\frac 12)^2, \beta=\gamma=1$ with the same initial condition as in \eqref{non:ini}.
\begin{comment}
is set to be
\begin{equation}\label{intial}
m(\bx,t)=(m_1,m_2,m_3)=
\begin{cases}
(0,0,-1)^T & \text{for} \quad   |\bx|\ge \frac 12,\\
(\frac{2x_1A}{A^2+|\bx|^2},\frac{2x_2A}{A^2+|\bx|^2},\frac{A^2-|\bx|^2}{A^2+|\bx|^2})^T & \text{for} \quad   |\bx| \leq \frac 12,
\end{cases}
\end{equation}
where $A=(1-2|\bx|)^4$. 
\end{comment}
We choose periodic boundary condition in space and implement Fourier spectral method for spatial discretization with Fourier modes $N=64$ in each direction.  Numerical simulations are computed by Gauss-Seidel scheme \eqref{LLG:correction:1b}-\eqref{LLG:correction:3b} with  $\delta t=10^{-5}$. 
Displayed in Fig.\;\ref{Fourier} are the orthogonal projection of the vector field $m$ on the plane $\{(x,y,0): x,y \in R\}$ at various times. We observe that $m$ changes rapidly near origin, indicating that   $\Grad m$  becomes very large and may blowup. It is observed from Fig.\;\ref{Fourier_3d} that the spin $m$ at the origin even change from $(0,0,1)$ to $(0,0,-1)$ at $t\approx 0.4$. %This indicate that our numerical schemes can reproduce the phenomenon of blowup for Landau-Lifshitz equation \eqref{Lf}. 
In Fig.\;\ref{grad_m}.(b),  we plot the  evolution of  $\|\Grad m\|_{\infty}$ with respect to time,  which is consistent with the phenomenon of blowup presented in \cite{bartels2008numerical}.

Using the same initial condition and resolution as above, we also plot  in Fig.\;\ref{grad_m}.(b) energy curves computed by the fourth-order Runge-Kutta scheme with $\delta t=10^{-5}$ and  the stabilized type-II Gauss-Seidel scheme with $S=0.5$ for the general Landau-Lifshitz equation using a larger time step $\delta t=10^{-2}$. %The type-II Gauss-Seidel scheme can be obtained by easily modified by type-I Gauss-Seidel scheme \eqref{LLG:correction:1b}-\eqref{LLG:correction:3b}.
We observe from  Fig.\;\ref{grad_m}.(a) that the two energy curves are close to each other  which indicates  that the stabilization allows us to use much larger time steps for the simulation of the general Landau-Lifshitz equation \eqref{Lf}.

\section{Concluding remarks}
We developed in this paper two classes of length preserving, and energy dissipative time discretization schemes for  the Landau-Lifshitz equation based on the Lagrange multiplier approach. These schemes are not restricted yo any particular spatial discretization, and at each time, the computational cost is dominated by the cost of the predictor step which requires solving  decoupled linear elliptic equations with constant coefficients. So the computational costs of these schemes  are comparable to a usual semi-implicit scheme, but with the advantage of being length preserving, and can also be energy dissipative.
 To the best of our knowledge, our schemes based on the predictor-corrector approach are the first  length preserving higher than second-order schemes for the Landau-Lifshitz equation, and our schemes with an additional space-independent Lagrange multiplier are the first  length preserving and energy dissipative schemes for the Landau-Lifshitz equation.
 
 We carried out ample numerical experiments  to validate the stability and accuracy for the proposed  schemes, and also compared them with some existing schemes. It is observed that 
 our predictor-corrector schemes  can provide better accuracy than  schemes based on projection. We also observed that the Type-II schemes, which are based on a  Lagrange multiplier $\lambda=0$ at the continuous level, 
usually   produce more accurate dynamical approximation than the type-I schemes, which are based on the Lagrange multiplier formulation with $\lambda=|\nabla m|^2$. %While type-I schemes may allow larger time steps than type-II schemes.  
It is also found that for the general Landau-Lifshitz equation, adding a stabilized term can significantly increase the stability for   both type-I and type-II schemes.

The general ideas introduced in this paper for constructing length preserving and energy dissipative schemes are not limited to Landau-Lifshitz equation, and can be applied to other length preserving nonlinear dissipative systems, such as the liquid crystal flows \cite{liu2000approximation,liu2002mixed}. In this sense, this paper can also be regarded as Part III of the sequence following \cite{CS_CAMME22,CS_SINUM22}.

We only provided some stability analysis for the simplest first-order scheme in the semi-discrete form. It appears very difficult to establish unconditional stability results for the higher-order schemes  in the semi-discrete form. In a future work, we  shall consider a fully discretized version of our schemes presented here, and attempt to establish its stability and convergence results with the help of some reasonable conditions on the time step, similar to those in \cite{an2021optimal,gui2022convergence}.

\bibliographystyle{siamplain}
\bibliography{references}

\end{document}